\documentclass[11pt,longtable,letterpaper,oneside,reqno]{amsart}
\usepackage{amsmath,amsfonts,amsthm,amssymb,amscd,mathtools,stmaryrd,url}
\usepackage{bbm}
\usepackage{epsfig}
\usepackage{graphicx}
\usepackage{latexsym}
\usepackage{longtable}
\usepackage{float}

\DeclareFontFamily{OT1}{rsfs}{}
\DeclareFontShape{OT1}{rsfs}{n}{it}{<-> rsfs10}{}
\DeclareMathAlphabet{\mathscr}{OT1}{rsfs}{n}{it}

\newcommand{\C}{\mathbb{C}}

\newcommand{\N}{{\mathbb{N}}}

\renewcommand{\L}{{\mathcal{L}}}
\renewcommand{\Re}{{\mathfrak{Re}}}
\renewcommand{\Im}{{\mathfrak{Im}}}

\renewcommand{\b}{\beta}
\renewcommand{\d}{\delta}
\newcommand{\e}{\eta}
\newcommand{\g}{\gamma}

\renewcommand{\k}{\kappa}
\renewcommand{\l}{\lambda}
\newcommand{\s}{\sigma}
\renewcommand{\t}{\theta}
\newcommand{\w}{\omega}
\newcommand{\ep}{\epsilon}

\newtheorem{prop}{Proposition}[section]
\newtheorem{thm}[prop]{Theorem}

\newtheorem{lem}[prop]{Lemma}

\newtheorem {defn }{Definition}
\newenvironment{rem}{{\bf Remark.}}{}
\numberwithin{equation}{section}
\usepackage[margin=1in]{geometry}
\usepackage{lipsum}
\begin{document}
\title[Explicit zero-free regions for Dirichlet $L$-functions]{Explicit zero-free regions for Dirichlet $L$-functions}
\author[H. Kadiri]{Habiba Kadiri}
\address{Department of Mathematics and Computer Science, University of Lethbridge, 4401 University Drive, Lethbridge, Alberta, T1K 3M4 Canada}
\email{habiba.kadiri@uleth.ca}
\thanks{
Research for this article is partially supported by an NSERC Discovery grant.
The calculations were executed on the University of Lethbridge Number Theory Group Eudoxus machine, supported by an NSERC RTI grant.}
\subjclass[2010]{11M26, 11M06}
\keywords{\noindent zero-free region, Dirichlet $L$-functions, explicit formula}
\maketitle
\dedicatory{{\it \`A mon p\`ere.}}
\begin{abstract}
Let $L(s,\chi)$ be the Dirichlet $L$-function associated to a non-principal primitive character $\chi$ modulo $q$ with $3\le q \le 400\,000$.
We prove a new explicit zero-free region for $L(s,\chi)$:
$L(s,\chi)$ does not vanish in the region $\Re s \ge 1-\frac1{R \log \left(q\max(1,| \Im s|)\right) }$ with $R=5.60$.
This improves a result of McCurley where $9.65$ was shown to be an admissible value for $R$.
\end{abstract}
\section{Introduction}
\label{intro}  
In this article, we establish an explicit zero-free region for the Dirichlet $L$-functions associated to moduli for which the Generalized Riemann Hypothesis has been partially verified.
\begin{thm}
\label{main-thm}
Let $q$ be an integer with $3\le q \le 400\,000$ and $\chi$ a non-principal primitive character modulo $q$. 
Then the Dirichlet $L$-function $L(s,\chi)$ does not vanish in the region: 
\begin{equation}\label{zfr-L}
\Re s \ge 1-\frac1{5.60 \log \left(q\max(1,| \Im s|)\right) }.  
\end{equation}
\end{thm}
This result improves previous results of \cite{Kad0} ($R=6.436$) and \cite{Mc2} ($R=9.646$).  
The case $q>400\,000$ will be treated in a follow-up article. 
Theorem \ref{main-thm} uses ideas developed in \cite{Kad1} where it is proven that the Riemann zeta function does not vanish in the region
\begin{equation}
\label{zfr-zeta} \Re s \ge 1-\frac{1}{5.70 \log |\Im s|},\  |\Im s|\ge2 .
\end{equation}
This improved long-standing results of Stechkin \cite{St} and Rosser and  Schoenfeld \cite{RS2} (they proved that a value just under $9.65$ was admissible).
Since then, the method of \cite{Kad1} to prove \eqref{zfr-zeta} has been refined and $5.70$ has been reduced to $5.68371$ in \cite{JanKwo} and to $5.573412$ in \cite{MosTru}.
\newline
Explicit zero-free regions for $L$-functions are useful 
to establish explicit results about the primes.
For instance, \eqref{zfr-zeta} was used to obtain estimates for finite sums or products over the primes, locating primes in short intervals, applications to Diophantine problems, etc
(see 
\cite{But1}, 
\cite{Dus2}, 
\cite{FaKa}, 
\cite{HajSarTij}, 
\cite{KaLu1},
\cite{Ram2}, 
and \cite{Tru4}). 
Recently, Theorem \ref{main-thm} has already been applied to obtain new explicit bounds for $\psi(x;q,a)$ in \cite{KaLu2} \cite{BeMaOBRe} (these works improve the well-known article \cite{RR} of Ramar\'e and Rumely).
Some of the ideas to establish \eqref{zfr-zeta} and \eqref{zfr-L} have already been extended to Dedekind zeta functions \cite{Kad3} and Hecke $L$-functions \cite{AhnKwo}. 
We now present a quick overview of the main ingredients of the proof.
\newline
Let $q$ be an integer with $3\le q \le 400\,000$ and $\chi$ a non-principal primitive character modulo $q$. 
We consider $\varrho_0 = \b_0 + i \g_0$ a non-trivial zero ($0<\b_0<1$) of the Dirichlet $L$-function $L(s,\chi)$.
Let $s$ be a complex number with $\Re s >0$, and let denote $\k$ and $\d$ some positive parameters.
Let $f$ be a compactly supported, non-negative, ``smooth'' function (it has continuous derivatives up to a certain order).
Instead of $-\Re \frac{L'}{L}(s,\chi) =  \Re \sum_{n\ge1} \frac{\Lambda(n)\chi(n) }{n^{s}} $, we consider the ``smoothed version'' 
\begin{equation}\label{eq0}
 \Re \sum_{n\ge1} \frac{\Lambda(n)\chi(n)f(\log n) \left(1- \frac{\k  }{n^{\d}} \right)}{n^{s}} .
\end{equation}
We establish a version of Guinand-Weil's explicit formula of the form
\begin{multline}\label{eq1}
 \Re \sum_{n\ge1} \frac{\Lambda(n)\chi(n)f(\log n) \left(1- \frac{\k  }{n^{\d}} \right)}{n^{s}}
\\=\frac{(1-\k  )}2  f(0) \log (q |\Im s| ) 
- \sum_{\varrho \in Z(\chi)} \Re \left(F(s-\varrho) -\k   F(s+\d-\varrho) \right)
 + E_q(s) ,
\end{multline}
where $F$ is the Laplace transform of $f$, 
$\chi$ is non-principal, 
$Z(\chi)$ is the set of non-trivial zeros of $L(s,\chi)$, and
$E_q(s)$ is an error term.
In addition, when $\chi$ is principal, the term $a_0\Re F(s-1) $ arises for $k=0$ from the pole of $\zeta(s) $ at $s=1$.
To compare with the classical proof, $\kappa$ and $\delta$ would each be $0$, $f$ would be identically $1$, $\Re F(s-1) $ would be $\frac1{\Re s-1}$, and $- \sum_{\varrho \in Z(\chi)} \Re\frac1{s-\varrho}$ would be  the sum over the zeros.
We apply \eqref{eq0} at various values of $s$ on a vertical line passing near $\varrho_0$.
In particular we take $s= \s+ik\g_0$ with the integer $k$ ranging between $0$ and a fixed integer $n_0$.
By means of a trigonometric inequality of the form
\[
P(t) = \sum_{k=0}^{n_0}  a_k \cos(kt) \ge 0\ \text{ with }\ a_k \ge0\ \text{ for all }\ k=0,\ldots, n_0,
\]
we deduce
\[
 \sum_{n\ge1} \frac{\Lambda(n)f(\log n) \left(1- \frac{\k  }{n^{\d}} \right)}{n^{\s}} \sum_{k=0}^{n_0}a_k \cos\left(k \arg\left(\frac{\chi(n)}{n^{i\g_0}} \right)\right) \ge 0.
\]
It remains to give accurate upper bounds to the right hand side of \eqref{eq1} for each $s=\s+ik\g_0$. 
First  we isolate the term arising from $k=1$ and the zero $\varrho_0$ in the sum over the zeros. 
In the classical proof, the sum over the remaining zeros is simply discarded as $\Re\frac1{s-\varrho}$ is always positive assuming $\Re s>1$. 
The situation here is more complicated as we are considering the difference $ \Re \left(F(s-\varrho) -\k   F(s+\d-\varrho) \right)$ and as 
we are allowing $s$ to be located inside the critical strip.
The main idea of Stechkin \cite{St} was to impose conditions of $\k,\d$ to deal with $ \Re \left(\frac1{s-\varrho} - \frac{\k}{s+\d-\varrho} \right)$.
This trick was then exploited by Rosser and Schoenfeld in \cite{RS2} and by McCurley in \cite{Mc2}.
We are using here the version adapted to smooth functions from \cite{Kad1}:
after imposing on $f,\k$ and $\d$ some appropriate conditions, we show
$
\displaystyle{ \sum_{\begin{substack}{\varrho \in Z(\chi) \\ \Re \varrho \le \Re s}\end{substack}} \Re \left(F(s-\varrho) -\k   F(s+\d-\varrho) \right)\ge 0 . }
$
We then prove that the size of the sum over the remaining zeros is negligible. 
This argument allows us to multiply the final constant in the zero-free region by a factor of $\frac{1-\k  }2 \simeq \frac{1-\frac1{\sqrt5}}2 \simeq0.28$. 
Note that this argument is valid for all moduli $ q$. 
On the other hand, the Burgess bound argument used by Heath-Brown in \cite{HB} leads to a coefficient of $0.25$ but was only valid for sufficiently large moduli.
Putting together all these arguments leads to the inequality:
\begin{equation}\label{unlab}
 \frac{1-\k  }2 f(0)\log\left(q |\g_0|\right)\sum_{k=1}^{n_0}a_k -  a_1 F(\s-\b_0) +a_0 F(\s-1)+ \epsilon  \ge 0,
\end{equation}
where $\epsilon$ is an error term.
We choose $f$ to depend on $\b_{0}$ by setting $f(0) =h(0) (1-\b_0)$, where $h(0)$ is independent of $\varrho_0$ and $h$ is a smooth function chosen appropriately. 
We also choose the polynomial coefficients $a_i$, and the parameter $\s$.
Then the inequality 
\begin{equation}\label{def-R0}
(1-\b_0)  \log\left(q  |\g_0| \right) \ge  \frac{ a_1  F(\s-\b_0) - a_0 F(\s-1)-\epsilon}{\frac{1-\k  }{2} h(0) \sum_{k=1}^{n_0}a_k}
\end{equation}
provides a formula where the zero-free constant $R^{-1}$ is given by the term on the right. 
\section{Explicit results about the zeros of Dirichlet $L$-functions}
We list here the most recent results for zeros of Dirichlet $L$-functions which shall be applied in this article.
\label{explicit-results-zeros-Dirichlet}
\subsection{Partial numerical verification of GRH}\label{Platt-RH}
In 2013, Platt provided a partial numerical verification of the Generalized Riemann Hypothesis which asserts that all zeros of Dirichlet $L$-functions inside the critical strip $0<\Re s <1$ lie on the vertical $1/2$-line:
\begin{thm}\cite[Theorem 7.1]{Pla1} \label{Platt}
GRH holds for Dirichlet $L$-functions of primitive characters modulo $3\le q\le 400\,000$ and to height 
$H_q=\max\left(\frac{10^8}q,\frac{c_q\cdot10^7}q+200\right) $
with $c_q = 7.5$ if $q$ is even and $3.75$ otherwise.
\end{thm}
This improves drastically both numerically and theoretically on Rumely's work (see \cite{RR} and \cite{Ben}).
It increases by a factor of $1\,000$ the number of moduli for which the verification was undertaken, and it increases by a factor of between $100$ and $10\,000$ the size of the height $H_q$ (for comparable moduli).
This theorem has many important applications and has already contributed to Helfgott's proof of the Ternary Goldbach Conjecture \cite{HH1} and \cite{HH2}.
\subsection{Explicit zero-free regions}
\begin{thm}\cite[Theorem 1]{Mc2}\label{zfr}
Let $\mathfrak{L}_q(s)$ be the product of the $\phi(q)$ Dirichlet $L$-functions formed with characters modulo $q$.
Let $M=\max\{q,q|\Im s|,10\}$ and $R=9.645908801$.  
Then $\mathfrak{L}_q(s)$ has at most a single zero in the region 
\[\left\{ s: \Re s \ge 1-\frac{1}{R \log  M } \right\}.\] 
The only possible zero in this region is a simple real zero arising from an $L$-function formed with a real non-principal character modulo $q$.
\end{thm}
In the PhD dissertation \cite{Kad0} it was established that $6.4355$ is an admissible constant for $R$.
\subsection{Explicit estimate of the number of zeros in a box}
We denote $N(T,\chi)$ be the number of zeros of $L(s,\chi)$ in the rectangle $0<\Re s<1$ and $|\Im s| < T$.
\begin{thm} 
\label{bound-N}
Let $\chi$ be a primitive non-principal character modulo $q$. Then there exist $C_1,C_2>0$ s.t. for every $T\ge 10$
\[
\left|N(T,\chi)-\frac{T}{\pi} \log \left(\frac{qT}{2\pi e}\right)\right| \le C_1 \log(qT) + C_2.
\]
In other words:
\[
N_2(T,q) \le N(T,\chi) \le N_1(T,q),\ 
\text{ with }\ 
N_j(T,q) = \frac{T}{\pi}\log \frac{qT}{2\pi e} + (-1)^{j+1}\left( C_1 \log (qT) + C_2 \right).
\]
\end{thm}
Theorem \ref{bound-N} was first established by McCurley \cite[Theorem 2.1]{Mc2}.
We use here Trudgian's version \cite[Theorem 1]{Tru1} with 
\begin{equation} \label{def-C1C2} C_1=0.247,\  C_2=8.949.\end{equation}
\section{Setting up the argument}
\label{set-up}  
\subsection{Notation}
\label{notations}
Let $q$ be a modulus for which the verification of the GRH up to height $H_q$ has been established.
In this article, we use Theorem \ref{Platt} and assume that 
\[2 \le q \le 400\, 000,\ H_q > 293,\ \text{ and }\  qH_q \ge Q_0 = 10^8.\]
Let $q_0, q_1$ be positive fixed integers so that 
$ q_0 \le q \le q_1.$
For computational purposes we split in two cases $3 \le q \le 1\,000$ and $1\,001 \le q \le 400\,000$, so
$q_0=3, q_1=1\,000, H_q\ge100\,000$ and $q_0=1\,000, q_1=400\,000, H_q\ge293$ respectively.
\newline
Let $\chi$ be a non-principal primitive character of conductor $q$
and let $\varrho_0 = \beta_0 + i\gamma_0$ be a non-trivial zero of $L(s,\chi)$ satisfying 
\begin{equation}\label{cond-gamma0}
\g_0 \ge H_q ,\ q \g_0 \ge Q_0, \ \hbox{ and }\ \b_0 < 1-\frac{1}{R \log  \left( q |\Im \g_0| \right)}\ \text{ with }\ R=9.646.
\end{equation}
We introduce the parameters $r$ and $\e$ such that 
\begin{equation}
 \label{def-eta}
5 \le r \le R \ \text{and}\ \e = 1-\beta_0= \frac{1}{r\log ( q  \g_0 )}.
 \end{equation}
Note that $r<5$ proves a zero-free region with admissible constant $5$.
\newline
Let $P$ be a trigonometric polynomial of degree $n_0\ge 2$ satisfying
\begin{equation}\label{cond-pol-trig} 
P(x) = \sum_{k=0}^{n_0} a_k \cos(kx) \ge0 , \ \hbox{ with }\ a_k\ge0, \ \text{and} \ 0<a_0<a_1.
\end{equation}
We denote
$
 A = \sum_{k=1}^{n_0} a_k.
$
We refer to Section \ref{choice-pol-trig} for the explicit definition of the polynomial used in this article. 
Here we use a polynomial of degree $n_0=16$.
\newline
Let $t_0>0$ (here $t_0=100$) and let 
\begin{align}
& \label{def-sigma}
 \s = 1-\frac1{R \log (q(n_0 \g_0 + t_0))},
  \\ & \label{def-omega}
  \omega=  \omega(\e,q) = \frac{1-\s}{\e} 
  = \frac{1}{R\e \log \left(q (n_0 \g_0 + t_0)\right)} 
  = \frac{1}{R\e \log \left(n_0 e^{\frac1{r\e}} + qt_0\right)} .
\end{align}
Since $q\g_0 \ge Q_0$, we have the ranges
\begin{equation}\label{def-sigma0}
 \max(\b_0,\s_0) < \s < 1\ \hbox{ with }  \s_0 =1 - \frac1{R\log\left(Q_0 n_0 + q_0t_0\right)},
 \end{equation}
 and
\begin{equation}\label{def-eta0}
 0< \e \le  \e_0 \ \hbox{ with }  \e_0 =\frac1{r\log  Q_0}.
\end{equation}
Let $\e_1\in(0,\e_0)$. We have
\begin{equation}
\label{ineq-omega} 
\begin{cases}
& \w_{1} \le \w \le \frac{r}R \ \hbox{ when }\ 0< \e \le \e_1,\\
& \w_{0} \le \w \le \w_1 \ \hbox{ when }\ \e_1< \e \le \e_0,
\end{cases}
\end{equation}
\begin{equation} \label{def-omega01}
\text{ with }\   \omega_0 =\omega(\e_0,q_1) = \frac{ r \log Q_0 }{ R \log\left(n_0 Q_0 + t_0 q_1\right)}
  \ \text{ and }\ 
  \omega_1 =\omega(\e_1,q_0)  = \frac{1}{R\e_1 \log \left(n_0 e^{\frac1{r\e_1}} + q_0t_0\right)}.
\end{equation}
Let $\d, \k  $ be some parameters satisfying
\begin{equation}
\label{def-kappa-delta}
0.5< \d<0.75,\ 
\text{and}\  0.25 < \k   < \min\left(\frac{\s_0}{1+\d} , \frac1{2\d+1}\right). 
\end{equation}
In addition, we assume
\begin{equation}\label{choice-t0}
\max \left( 10, \left( \left( \k    (2\d +1) \right)^{-1}-1\right)^{-1/2}\right) < t_0 < 293< H_q.
\end{equation}
We provide in Proposition \ref{gal-Stechkin} the definitions of $\k  $ and $\d$ depending on $h$, $\s_0$, and $\eta_0$.
\newline
Finally we add the condition
\begin{equation}\label{cond-primitivity-kappa-delta2}
a_0 +  \left( \frac{1-\kappa}{\max\left(
\frac1{2^{\s_0}-1}-\frac{\k}{2^{\s_0+\d}-1}
,2
\left(\frac1{3^{\s_0}-1}-\frac{\k}{3^{\s_0+\d}-1}\right)
\right)} - 1 \right) \sum_{k=2}^{n_0}  a_k  >0 .
\end{equation}
Numerical data for all the above parameters ($r,\k,\d $) are given in Tables \ref{table-R01} and \ref{table-R02}, Section \ref{Computations}.
\subsection{Introducing a smooth weight.}
\label{smooth-weight}
We now introduce the weight $f$ which is used in the study of $\sum_{n\ge1}\frac{\chi(n)  \Lambda(n)}{n^s}f(\log n)\left(1-\frac{\k}{n^{\d}}\right) $.
We define
\begin{equation}
f(t) = f_{\e}(t) = \e h(\e t),
\end{equation}
where $h$ is a function satisfying:
\begin{equation}\label{cond-f} 
\left\{\begin{split}
& h \text{ is compactly supported in $[0,d]$, for some }\ d>0,\\
& h \text{ is positive in } [0,d),\\
& h \in \mathcal{C}^2([0,d]), \\
& h(d) = h'(0) = h'(d) = h^{\prime \prime}(d) = 0.
\end{split}\right.
\end{equation}
Note that in \cite{Kad1}, we used the notation $g_1(\t)$ and $m$ instead of $h(0)$ and $m_h $ respectively.
We give more details in Section \ref{choice-weight} on the motivation for the explicit choice for $f$ as well as study its analytical properties.
\newline
We denote $F$ the Laplace transform of $f$:
 \begin{equation}\label{def-Laplace-F}
 F(s) 
= \int_{0}^{\infty} e^{-st} f(t) dt
=  \int_{0}^{d } e^{-\frac{s u}{\e}}   h(u) d u,
 \end{equation}
$F_2$ the Laplace transform of $f^{\prime\prime}$, and
\begin{equation}\label{def-H}
H(x,y)= \Re \left( \frac{F_2(x+iy)}{(x+iy)^2}\right)  .
 \end{equation}
We impose a last and important condition on $f$: that its Laplace transform satisfies
\begin{equation}
 \label{pos-F}
 \Re F (z) \ge 0 \ \text{ when  }\ \Re z \ge 0.
\end{equation}
\begin{rem}
This is achieved by choosing $h$ as a self convolution: $h=g\star g$, and we refer the reader to \cite[Section 7]{HB} for details.
In the classical proof, property \eqref{pos-F} is satisfied by $\frac1s$ and is key in handling the sum over the zeros. 
\end{rem}
\newline
We now recall some properties established in \cite[Lemma 3.2]{Kad1}
\begin{lem}\label{lem-estimates-F}
Let $x,y$ be real numbers satisfying $x+iy\not=0$.
\begin{equation}\label{dvpt-F}
\Re F(x+iy) = \frac{x}{x^2+y^2} h(0) \e + H(x,y),
\end{equation}
 \begin{equation}\label{boundF-H}
\text{ where }\ |H(x,y)| \le \frac{M\left(x/\e \right)}{x^2+y^2}  \e^2, \text{ with }  M(z) = \int_{0}^{d } \left|h''(u)\right| e^{-zu} du .
\end{equation}
In addition, when $x \ge 0$, then
\begin{equation}\label{boundF-m}
|H(x,y)| \le  \frac{m_h }{x(x^2+y^2)}  \e^3, \text{ with }   m_h  = \max _{u\in[0,d ]} |h^{\prime\prime}(u)|.
\end{equation}
\end{lem}
\begin{rem}
We choose $h$ among several choices of families of functions. The explicit definition of $h$ eventually depends on several  extra parameters (see Section \ref{choice-weight}).
\end{rem}
\subsection{An explicit formula for a smoothed version of $-\frac{L'}{L}(s,\chi)$.}
\label{Explicit-Formula-Lfn}
\begin{prop}\label{Prop1.1}
Let $f$ be a function satisfying \eqref{cond-f}. Let $s$ be a complex number, and $\chi$ a primitive (non-principal) character modulo $q$. Then
\begin{multline} \label{Explicit-formula-primitive-f}
\Re \sum_{n\ge1}\frac{\Lambda(n)\chi(n) f(\log n) }{n^s} =
- \Re \sum_{\varrho \in Z(\chi)} F(s-\varrho)
+\frac{ f(0)}2\left( \log \frac{q}{\pi} + \Re \frac{\Gamma'}{\Gamma} \left(\frac{s+\mathfrak{a}}{2}\right) \right)
\hfill \\
+\Re \frac1{2i\pi} \int_{\frac12-i\infty}^{\frac12+i\infty} \Re \frac{\Gamma'}{\Gamma} \left(\frac{z+\mathfrak{a}}{2}\right)\frac{F_2(s-z)}{(s-z)^2} d z,
\end{multline} 
\begin{multline} \label{Explicit-formula-zeta-f}
\Re \sum_{n\ge1}\frac{\Lambda(n) f(\log n)  }{n^s}
=
 \Re F(s-1) 
 - \sum_{\varrho \in Z(\zeta)} \Re F(s-\varrho) 
+\frac{ f(0)}2 \left(  -  \log \pi + \Re \frac{\Gamma'}{\Gamma} \left(\frac{s}{2}+1\right)  \right)
\hfill \\  + \Re \frac{F_2(s)}{s^2}
+ \Re \frac1{2i\pi} \int_{\frac12-i\infty}^{\frac12+i\infty} \Re \frac{\Gamma'}{\Gamma} \left(\frac{z}{2}\right)\frac{F_2(s-z)}{(s-z)^2}d z,
\end{multline}
where $Z(\zeta)$ and $Z(\chi)$ are the sets of non-trivial zeros of $\zeta(s)$ and of $L(s,\chi)$, respectively. 
\end{prop}
Proof of Proposition \ref{Prop1.1} is postponed until Section \ref{proof-ExplicitFormula}.
The case of the principal character $\chi_0$ follows from \eqref{Explicit-formula-zeta-f} and the identity
\[
\Re \sum_{n\ge1}\frac{\Lambda(n)\chi_0(n) f(\log n)}{n^s}
=\Re
\sum_{n\ge1}\frac{\Lambda(n)f(\log n) }{n^s}
- \Re \sum_{m\ge1} \sum_{p\mid q}\frac{(\log p) f(m\log p)}{p^{ms}} .
\]
\newline
Let $s=x+iy$ be a complex number. 
For $\k,\d$ as defined in \eqref{def-kappa-delta}, we introduce
\begin{align}
& \label{def-S} S(s) = \Re \sum_{n\ge1} \frac{\Lambda(n) f(\log n)  }{n^s}\left(1- \frac{\k  }{n^{\d}}\right) ,
\\& \label{def-Schi} S(s,\chi) = \Re \sum_{n\ge1} \frac{\chi(n)\Lambda(n)f(\log n)  }{n^s}\left(1- \frac{\k  }{n^{\d}}\right) .
\\&
\label{def-E}
E(x,y) = H(x,y) - \k   H(x+\d,y) .
\\& 
\label{def-D}
D(s) =\Re(F (s) - \k  F(s+\d) )=  \Big( \frac{x}{x^2+y^2} -  \frac{\k(x+\d)}{(x+\d)^2+y^2}\Big)h(0) \e +  E(x,y),
\\&
\label{def-real-log-der-Gamma}
\mathfrak{G}_{\k,\d}(s) =  \Re \frac{\Gamma'}{\Gamma}(s/2) - \k  \Re\frac{\Gamma'}{\Gamma}((s+\d)/2) .
\end{align}
In addition we define
\begin{align} 
\label{def-D1}
& D_1(s)  = \frac{h(0) \e}2 \left[ -(1-\k) \log \pi  
+\mathfrak{G}_{\k,\d} (s+2 ) 
\right],
\\&
\label{def-D1chi}
D_1(s,\chi) 
 = \frac{h(0) \e}2 \left[ (1-\k) \log{\frac{q}{\pi}} 
+\mathfrak{G}_{\k,\d} ( s+\mathfrak{a}) 
\right],
\\&
\label{def-D2} 
D_2(s) =   E(x,y)
+ \frac1{2\pi} \int_{-\infty}^{+\infty}  \Re \frac{\Gamma'}{\Gamma}( (1/2 +it)/2 ) E(x-1/2 ,y-t) dt , 
\\&
\label{def-D2chi}
D_2(s,\chi) = \frac1{2\pi} \int_{-\infty}^{+\infty} \Re \frac{\Gamma'}{\Gamma} ( (1/2 +\mathfrak{a}+it)/2 )  E(x-1/2 ,y-t) dt .
\end{align}
Using this notation and noting that  $f(0)=h(0)\e$ and $\Re \frac{\Gamma'}{\Gamma}(s) =  \Re \frac{\Gamma'}{\Gamma}(s/2)$, the explicit formulae from Proposition \ref{Prop1.1} become
\begin{align}
\label{notation-expl-form-S}
& S(s) = D(s-1) - \sum_{\varrho \in Z(\zeta)} D(s-\varrho)  +  D_1(s) + D_2(s),
\\
\label{notation-expl-form-Schi}
& S(s,\chi) = - \sum_{\varrho \in Z(\chi)} D(s-\varrho) + D_1(s,\chi) + D_2(s,\chi) .
\end{align}
\subsection{An explicit inequality}\label{trig-pol}
Taking $x=\arg\left(\frac{\chi (n)}{n^{i \g_0}}\right)$ in the non-negative trigonometric inequality \eqref{cond-pol-trig}, we get
\begin{equation}\label{Pos1}
\sum_{k=0}^{n_0} a_k S(\s+ik\gamma_0,\chi^k)
 \ge 0.  
\end{equation}
Note that Proposition \ref{Prop1.1} can be applied to $\chi^k$ only when it is a non-primitive character. 
So for each $k\ge2$, we introduce $q_k$ the conductor for the Dirichlet character $\chi^k$, and 
$\chi_{(k)}$ the unique primitive character modulo $q_k$ that induces $\chi^k$. 
(We refer the reader to \cite[Theorem 9.2.]{Mont} and recall that for $k=0$, $\chi_{(0)}$ is identically $1$, and for $k=1, \chi_{(1)}=\chi$.)
We establish in Lemma \ref{primitivity} that under condition \eqref{cond-primitivity-kappa-delta2}, then 
\begin{equation}\label{Pos2}
\sum_{k=0}^{n_0}  a_k S(\sigma+ik\gamma_0,\chi_{(k)}-\chi^k)
+ \frac{1-\kappa}2 f(0) \sum_{k=1}^{n_0}  a_k \log \left(\frac{q}{q_k}\right)
\ge 0.
\end{equation}
Together with \eqref{Pos1} and \eqref{Pos2}, we obtain
\begin{equation}\label{pos-trig}  
\sum_{k=0}^{n_0}  a_k  S(\s+ik\gamma_0,\chi_{(k)}) + \frac{1-\k  }{2}f(0) \sum_{k=2}^{n_0}  a_k \log  \left(\frac{q}{q_k}\right)
\ge 0,
\end{equation}
and we can now apply the explicit formulae from Proposition \ref{Prop1.1} to each ``primitive term'' $S(\s+ik\gamma_0,\chi_{(k)}) $.
We introduce the notation
\begin{align}
\label{def-mathcalE1}
\mathcal{E}_1
& = a_0 D(\s-1) = a_0\left( F(\sigma-1,0)- \k   F(\sigma-1+ \d ,0)  \right), \hfill
\\ \label{def-mathcalE2}
\mathcal{E}_2  & = \sum_{k=0}^4 a_k \sum_{\varrho \in Z(\chi_{(k)})} D(\s+ik\g_0-\varrho), 
\\ \label{def-mathcalE3}
\mathcal{E}_3 
& =a_0 D_1(\s)+ \sum_{k=1}^{n_0}  a_k \left( D_1(\s+ik\g_0,\chi_{(k)}) + \frac{1-\k  }2f(0) \log ( q/ q_k)\right)
\\\nonumber & = 
\frac{h(0) \e}2  \Big[ 
a_0   \left(- (1-\k )\log \pi 
+  \mathfrak{G}_{\k,\d} ( \s+2 )
 \right)
 +  (1-\k) A  \log(q / \pi)    
+  \sum_{k=1}^{n_0}  a_k \mathfrak{G}_{\k,\d} (\s+\mathfrak{a}_k+ik\g_0) 
\Big],
\\ \label{def-mathcalE4}
\mathcal{E}_4 
&= a_0 D_2(\s) + \sum_{k=1}^{n_0}  a_k D_2(\s+ik\gamma_0,\chi_{(k)})  
\\\nonumber & = a_0  E(\s,0)
+ \sum_{k=0}^{n_0}  \frac{a_k}{2\pi}  \int_{-\infty}^{+\infty} \Re \frac{\Gamma'}{\Gamma} \left( \frac{1/2 +\mathfrak{a}_{k}+it}2 \right) E  (\s-1/2 ,k\g_0-t) 
 dt,
\end{align}
where
$\mathfrak{a}_k$ is associated to the character $\chi_k$: 
\begin{equation}\label{mathfrak-ak}
\mathfrak{a}_k=\frac{1-\chi_k(-1)}2. 
\end{equation}
In particular $\mathfrak{a}_0=0$.
Using this notation to rewrite explicit formulae \eqref{notation-expl-form-S} and \eqref{notation-expl-form-Schi} respectively, we rewrite \eqref{pos-trig} more simply as
\begin{equation}\label{pos-mathcalEs}
0\le  \mathcal{E}_1 -\mathcal{E}_2+ \mathcal{E}_3 + \mathcal{E}_4 .
\end{equation}
\section{Sketch of the proof of the Zero-Free Region}
\subsection{Studying each $\mathcal{E}_i$}
We make use of Lemma \ref{lem-estimates-F} for this and refer the reader to the next section for full details. 
\begin{itemize}
\item 
For $\mathcal{E}_1$, the pole of the Riemann zeta function at $s=1$ gives the contribution $a_0 D(\s-1) =a_0 F(\s-1) +\mathcal{O}(\e)= \frac{a_0}{\s-1}+\mathcal{O}(\e)$. 
We prove in Lemma \ref{lem-mathcalE1} that
\[
\mathcal{E}_1
\le  a_0 F(\s-1) +\mathcal{C}_{1}(\e).
\]
\item 
In the sum over the zeros $\mathcal{E}_2 $, we isolate the term at $k=1,\varrho=\b_0+i\g_0$:
\[
a_1\left(F (\s-\b_{0}) -\k F (\s+\d-\b_{0})\right) = a_1 F (\s-\b_{0}) + \mathcal{O}(\e).
\]
The positivity property \eqref{pos-F} for the Laplace transform $F$ allows us to discard most of the zeros on the right of the vertical line at $\s$.
This is done by means of a generalization of Stechkin's lemma (see Section \ref{section-kappa} and in particular Proposition \ref{gal-Stechkin}): \\
for $d$ the solution of the equation
\begin{multline}
\label{def-delta} 
(2t+1)h(0)+\left(\frac1{t}+\frac1{2\sigma_0-1+t}\right)m_h \e_0^2 
\\= \left(\frac1{t}+\frac{1+t}{(2\sigma_0-1+t)^2}\right)h(0) + \left(\frac1{t^3}+\frac1{(2\sigma_0-1+t)^3}\right) m_h \e_0^2 
\end{multline}
and for 
\begin{equation}\label{def-kappa}
\k=\frac{(2\sigma_0-1)h(0)-\frac{m_h \e_0^2}{2\sigma_0-1}}{(2\d +1)h(0) +\left(\frac1{\d }+\frac1{2\sigma_0-1+\d }\right)m_h \e_0^2},
\end{equation}
then \[
D(\sigma -\beta+i y)+D(\sigma-1+\beta+i y) \ge 0
\]
as soon as $1/2<\beta<\s$ and $y>0$.\\
Finally, we use an explicit zero density estimate to prove that the remaining  sum contributes negligibly. 
We prove in Lemma \ref{lem-mathcalE2} that
\begin{equation}
\mathcal{E}_2 
\ge a_1 F (\s-\b_{0},0) + \mathcal{C}_{2}(\e).
\end{equation}
 \item 
We use Stirling's formula to estimate the Gamma terms appearing in $\mathcal{E}_3$ and $\mathcal{E}_4$.
Since $\g_0>1$, we have 
$
\left|\mathfrak{G}_{\k,\d}(\s+\mathfrak{a}_k+ik\g_0) \right| \le \log \g_0  +\mathcal{O}(1)
$, 
for all $k\ge 1$.
We prove in Lemma \ref{lem-mathcalE3} that 
\[
\mathcal{E}_3
\le  \frac{A (1-\k)h(0)}2 \log (q \g_0) \e + \mathcal{C}_{3}(\e).
\]
 \item Finally, the integral term is negligible and we prove in Lemma \ref{lem-mathcalE4} that 
\[
\mathcal{E}_4 \le \mathcal{C}_{4}(\e).
\]
\end{itemize}
The positivity argument \eqref{pos-mathcalEs} can then be rewritten as
\begin{equation}\label{idea-final-pos}
 0 \le a_0 F(\s-1) - a_1 F (\s-\b_{0})+  \frac{A (1-\k)h(0)}2 \log (q \g_0) \e +\mathfrak{e}(\e),
\end{equation}
where the error term is 
\begin{equation}\label{def-error}
\mathfrak{e}(\e) = \mathcal{C}_{1}(\e)+\mathcal{C}_{2}(\e)+\mathcal{C}_{3}(\e)+\mathcal{C}_{4}(\e).
\end{equation}
We introduce the difference
\begin{equation}\label{def-intK}
K_h(\w) =  a_1   F(\s-\b_0) - a_0 F (\s-1)= \int_0^{d } \left(a_1 e^{-t} - a_0 \right) h(t) e^{\omega t} dt,
\end{equation}
so that inequality \eqref{idea-final-pos} becomes
\begin{equation}\label{ineq-R}
r = \frac1{\log (q \g_0) \e} \le \frac{A (1-\k)h(0)}{2(K_h(\w) - \mathfrak{e}(\e))}.
\end{equation}
\subsection{A strategy to compute the constant in the zero-free-region}
It remains to find a lower bound as large as possible for $K_h(\w) - \mathfrak{e}(\e)$ for the values of $\e$ in $(0,\e_0)$.
\newline
Since 
$
\displaystyle{ \frac{\partial K_h}{\partial \w}(\w) 
= \int_0^{d } \left(a_1 e^{-t} - a_0 \right) h(t) t e^{\omega t} dt,
}$
with $a_1 e^{-t} - a_0 >0$ when $t<\log\frac{a_1}{a_0}$ and $\w_{0} \le \w \le \frac{r}{R}$, then
$
\displaystyle{ 
\frac{\partial K_h}{\partial \w}(\w) 
> \int_0^{\log(a_1 / a_0)} \left(a_1 e^{-t} - a_0 \right) h(t)t e^{\w_{0} t} dt 
 + \int_{\log(a_1 / a_0)}^{d } \left(a_1 e^{-t} - a_0 \right) h(t)t e^{\frac{r}{R}t} dt .
}$
We choose values for $\w_{0}$ so that 
\begin{equation}
\label{pos-der-K}
\int_0^{\log(a_1 / a_0)} \left(a_1 e^{-t} - a_0 \right) h(t)t e^{\w_{0} t} dt + \int_{\log(a_1 / a_0)}^{d } \left(a_1 e^{-t} - a_0 \right) h(t)t e^{\frac{r}{R}t} dt
>0
\end{equation}
and $K_h(\w)$ increases with $\w \in (\w_0,r/R)$.
\newline
From the definitions of the $\mathcal{C}_i$'s in \eqref{newdef-C1}, \eqref{newdef-C2}, \eqref{newdef-C3}, and \eqref{newdef-C4}, we have
\begin{equation}
\mathfrak{e}(\e) =\e(\alpha_1 + \alpha_2 \e + \alpha_3 \e^2),
\end{equation}
where the $\alpha_i$'s are computable constants and satisfy $\alpha_1 <0, \alpha_2>0, \alpha_3>0$.
We have that $\mathfrak{e}(\e)$ is negative and decreases from $\e=0$ to $\eta_2 = \frac{-\alpha_2+\sqrt{\alpha_2^2-3\alpha_1\alpha_3}}{3\alpha_3}$, then it increases (becoming positive after the root $\frac{-\alpha_2+\sqrt{\alpha_2^2-4\alpha_1\alpha_3}}{2\alpha_3}$).
\newline
We are now able to establish a lower bound for $K_h(\w) - \mathfrak{e}(\e)$.
We fix $\eta_1\in(0,\eta_0)$ and bound $K_h(\w) - \mathfrak{e}(\e)$ in each following cases, depending on the location of $\eta_1$ with respect to $\eta_2$. 
We denote
\begin{equation}\label{def-omega12}
 \w_1 = \omega(\eta_1,q_0),\ \w_2 = \omega(\eta_2,q_0).
\end{equation}
\begin{enumerate}
 \item If $0<\eta_0<\eta_2$, 
then
\[ K_h(\w) - \mathfrak{e}(\e) > 
\begin{cases}
K_h(\w_1) &\ \mathrm{ if }\ 0<\eta<\eta_1,\\
K_h(\w_0) - \mathfrak{e}(\e_1)  &\ \mathrm{ if }\ \eta_1<\eta<\eta_0.
\end{cases}
\]
 \item If $0<\eta_2<\eta_0$ and $0<\eta_1<\eta_2$
then
\[ K_h(\w) - \mathfrak{e}(\e) > 
\begin{cases}
K_h(\w_1) &\ \mathrm{ if }\ 0<\eta<\eta_1,\\
K_h(\w_0) - \max(\mathfrak{e}(\e_1),\mathfrak{e}(\e_0))  &\ \mathrm{ if }\ \eta_1<\eta<\eta_0.
\end{cases}
\]
 \item If $0<\eta_2<\eta_0$ and $\eta_2<\eta_1<\eta_0$
then
\[ K_h(\w) - \mathfrak{e}(\e) > 
\begin{cases}
K_h(\w_1) &\ \mathrm{ if }\ 0<\eta<\eta_1,\\
K_h(\w_0) - \mathfrak{e}(\e_0)  &\ \mathrm{ if }\ \eta_1<\eta<\eta_0.
\end{cases}
\]
\end{enumerate}
We deduce that, for all $\eta\in (0,\e_0)$, $K_h(\w) - \mathfrak{e}(\e) >\mathcal{K}_1$, with
\begin{equation}\label{def-K1}
\text{  }\ 
\mathcal{K}_1 = \begin{cases}
\min\left( K_h(\w_1) ,K_h(\w_0) - \mathfrak{e}(\e_1)  \right) &\ \mathrm{ if }\ 0<\eta_0<\eta_2,\\
\min \left( K_h(\w_1) , K_h(\w_0) - \max(\mathfrak{e}(\e_1),\mathfrak{e}(\e_0)) \right)  &\ \mathrm{ if }\ 0<\eta_2<\eta_0.
\end{cases}
\end{equation}
Finally, we choose the value for $\eta_1<\eta_0$ so as to make $ \mathcal{K}_1$ as large as possible within the constraints of \eqref{def-eta}, \eqref{cond-pol-trig}, \eqref{def-kappa-delta}, \eqref{choice-t0}, \eqref{cond-primitivity-kappa-delta2}, and \eqref{cond-f}.
With this value for $ \mathcal{K}_1$, we define the constant in the zero-free region by
\begin{equation}\label{formula-R0}
\dfrac{A(1-\k  )h(0)}{ 2  \mathcal{K}_1}.
\end{equation}
Details for the computations can be found in Section \ref{Computations}.
\newline
\begin{rem}
In \cite{Kad1}, values for $\eta_0$ were chosen so that $\mathfrak{e(}\eta_0)<0$.
This allowed us to take 
\[\frac{A(1-\k  )h(0)}{ 2 K_h(\w_0)}\]
as an admissible value for the constant in the zero-free region for the Riemann zeta function.
The main improvement in \cite{MosTru} comes from refining this argument by dividing the interval of study at a value $\eta_1$ between $0$ and $\e_0$. 
The constant is then given by
\[\frac{A(1-\k  )h(0)}{ 2 }\min_{0< \e_1 \le \e_0} \max\left(\frac1{K_h(\w_1)},\frac1{K_h(\w_0)-\mathfrak{e(}\e_1)}\right),\]
assuming $\mathfrak{e(}\e)$ is decreasing on $(0,\e_0)$ and $\e_0\le \eta_2$. 
\end{rem}
\section{Details of the proof of the zero-free region}\label{details-proof}
\subsection{Study of polar  term $\mathcal{E}_1$ 
}\label{Study-mathcalE1}
\begin{lem}\label{lem-mathcalE1}
Let $\k  ,\d>0$, and $0<\s<1$. Then
\begin{equation}\label{bound-a0D}
\mathcal{E}_1 \le a_0 F (\s-1,0) + \mathcal{C}_{1} (\e),
\end{equation}
\begin{equation}\label{newdef-C1}
\text{ where }\ \mathcal{C}_{1} (\e) = a_0\left(- \frac{h(0) \k  }{\d} \e + \frac{m_h \k  }{(\s_0-1+\d)^3} \e^3\right).
\end{equation}
\end{lem}
\begin{proof}
We use \eqref{dvpt-F} and \eqref{boundF-m} to bound $F(\s-1+\d)$:
\[
\left| F(\s-1+\d,0) - \frac{h(0)}{\s-1+\d}\e \right|\le  \left| H(\s-1+\d,0) \right| \le \frac{m_h }{(\s-1+\d)^3}  \e^3 .
\]
\end{proof}
\subsection{Study of the sum over the zeros $\mathcal{E}_2 $} 
\label{Study-mathcalE2-sum-zeros} 
Let $\chi$ be a primitive character modulo  $q$ and let $k=0,...,n$.
We recall that $ \chi_{(k)}$ is the primitive character induced by $\chi^k$,
$q_k$ is its conductor, $Z(\chi_{(k)})$ is the set of non-trivial zeros of $L(s,\chi_{(k)})$, and in particular
$Z(\chi_{(0)})$ is the set of non-trivial zeros of the Riemann zeta function $\zeta(s)$. 
We study here
\begin{equation}\label{def-bigS}
\mathcal{E}_2  = \sum_{k=0}^{n_0}  a_k  
\sum_{\varrho \in Z(\chi_{(k)})} D(\s+ik\g_0-\varrho).
\end{equation}
Using the symmetry of the zeros with respect to the critical line, we write 
\begin{equation}\label{sym-zeros}
\sum_{\varrho \in Z(\chi_{(k)})} D(\s+ik\g_0-\varrho)
 = \sum_{\varrho  \in Z(\chi_{(k)}) }^{\star} \left[ D(\sigma -\beta+i(k\gamma_0-\gamma))+D
(\sigma-1+\beta+i(k\gamma_0-\gamma)) \right], 
\end{equation}
where $\sum_{\varrho}^{\star} = \frac12\sum_{\beta =\frac12 } + \sum_{\beta > \frac12 }$.
\subsubsection{Isolating the zero $\b_0+i\g_0$}
We isolate the summand for $k=1$ and $\varrho=\varrho_{0}$ 
and use the following inequality \cite[(34) page 325]{Kad1} :
\begin{equation}
D(\sigma-\beta_0) + D(\sigma-1+\beta_0) 
\ge F(\sigma-\beta_0,0) -  c_{2,1}(\e), 
\end{equation}
\begin{equation}\label{newdef-c21}
\text{ where }\  c_{2,1}(\e) = 
-\left[ 1-\k  \left(\frac1{\d}+\frac1{\s_0-\e_0+\d}\right)\right]h(0) \e
 + \left[ 1 + \k  \left(\frac1{\d^3}+\frac1{(\s_0-\e_0+\d)^3}\right)\right]m_h  \e^3.
\end{equation}
Thus\begin{multline}\label{ineq1-bigS}
\mathcal{E}_2 
\ge 
a_1  \left( F(\sigma-\beta_0,0) - c_{2,1}(\e) \right)
\\+ \Big( a_1  \sum_{ \begin{substack} {k=1\\\varrho  \in Z(\chi) \\ \varrho\not=\b_0+i\g_0 } \end{substack} }^{\star} 
+ \sum_{\begin{substack}{0\le k \le n\\k\not=1 }\end{substack}}  a_k   \sum_{\varrho  \in Z(\chi_{(k)}) }^{\star} \Big)
\left[ D(\sigma -\beta+i(k\gamma_0-\gamma))+D
(\sigma-1+\beta+i(k\gamma_0-\gamma)) \right].
\end{multline} 
One of the key arguments in this proof consists in reducing the size of the above sums.
Note that $D(\sigma -\beta+iy)$ has the same sign as essentially 
$\displaystyle{
\frac{\sigma -\beta}{(\sigma -\beta)^2+y^2} - \k \frac{\sigma -\beta+\delta}{(\sigma -\beta+\delta)^2+y^2} .
}$\\
In \cite[Lemma 2]{St} Stechkin prove that this is positive when $\sigma -\beta>0$ and under certain conditions for $\k,\d$.
This is done in \cite{Kad1} and we recall this result here.
\subsubsection{
A positivity argument for the zeros on the left of $\s$} \label{section-kappa}
\begin{prop}\cite[Proposition 4.2]{Kad1}\label{gal-Stechkin}
Let $1/2<\s_0<\s<1$ and $0<\e<\e_0$. Let $h$ be a positive function satisfying \eqref{cond-f}.
We define
\begin{align}
&\label{def-k2}
\k_2(t) =
\frac{(2\sigma_0-1)h(0)-\frac{m_h \e_0^2}{2\sigma_0-1}}
{(2t+1)h(0)
+\left(\frac1{t}+\frac1{2\sigma_0-1+t}\right)m_h \e_0^2},
\\ &\label{def-k3}
\k_3(t) =
\frac{(2\sigma_0-1)h(0)-\frac{m_h \e_0^2}{2\sigma_0-1}}{
\left(\frac1{t}+\frac{1+t}{(2\sigma_0-1+t)^2}\right)h(0)
+m_h \e_0^2 \left(\frac1{t^3}+\frac1{(2\sigma_0-1+t)^3}\right)}.
\end{align}
If $1/2<\beta<\s$ and $y>0$, then
$\displaystyle{
D(\sigma -\beta+i y)+D(\sigma-1+\beta+i y) \ge 0, 
}$ 
as long as $0\le x \le \min\left(\k_2(t),\k_3(t)\right)$ and
$t \ge \d$, where $\d$ is the solution in the interval $(0.5,0.75)$
of the equation 
$\k_2(t) = \k_3(t).$
We denote $\k $ the corresponding value of $\k_2$ at
$\d$: $\k =\kappa_2(\d)=\kappa_3(\d)$.
\end{prop}
Note that as $\s_0$ and $\e_0$ depend on $r$, it follows that $\k$ and $\d$ depend on $h$ and $r$.
\\
This proposition is key in the proof: it allows to reduce the final constant in the zero-free region by a factor of $(1-\k)$.
As a consequence of Proposition \ref{gal-Stechkin}, we can discard all zeros to the right of $\s$, so that
\begin{multline}\label{ineq21-bigS}
\mathcal{E}_2 
\ge a_1  \left( F(\sigma-\beta_0,0) - c_{2,1}(\e) \right)
\\ + \sum_{k=0}^{n_0}   a_k \sum_{{\begin{substack} {\varrho \in Z(\chi_{(k)}) \\ \b \ge \s } \end{substack}}}^{\star}\left[ D(\sigma -\beta+i(k\gamma_0-\gamma))+D(\sigma-1+\beta+i(k\gamma_0-\gamma)) \right].
\end{multline}
\subsubsection{Estimating the contribution of the zeros on the right of $\s$}
\label{Estimating-contribution-zeros-right-sigma}
In previous classical proofs of zero-free regions, the parameter $\s$ was greater than 1 as it was in the region of convergence of $\zeta(s)$. 
Here, the explicit formula for the smoothed version of $\Re\frac{L'}{L}$ allows us to choose our parameter $\s+i\g_0$ inside the critical strip, and thus closer to the zero $\varrho_0$ we need to locate.
This appears in $K(\w) = \int_0^{d } \left(a_1 e^{-t} - a_0 \right) h(t) e^{\w t} dt $ which needs to be as large as possible to reduce our final constant in the zero-free region.
This is feasible as we allow $\w = \frac{1-\s}{1-\b_0}$ to be positive.
On the other hand, a new contribution from the zeros in the vertical strip between $\s$ and $1$ arises in the sum over the zeros. 
We prove here that this one is indeed negligible.
From the zero-free region of Theorem \ref{zfr}, $\b \ge \s$ implies $|\g| \ge k\g_{0}+t_{0} $. Denoting
\begin{equation} \label{def1-bigSk}
\mathcal{E}_{2,k} 
= \sum_{{\begin{substack} {\varrho \in Z(\chi_{(k)}) \\ |\g| \ge k\g_0+t_0} \end{substack}}}^{\star}
\left[ D(\sigma -\beta+i(k\gamma_0-\gamma))+D(\sigma-1+\beta+i(k\gamma_0-\gamma)) \right] ,
\end{equation}
the inequality \eqref{ineq21-bigS} becomes 
\begin{equation}\label{ineq2-bigS}
\mathcal{E}_2 
\ge a_1  \left( F(\sigma-\beta_0,0) - c_{2,1}(\e) \right)
+ \sum_{k=0}^{n_0}   a_k \mathcal{E}_{2,k}.
\end{equation}
We now use the estimate \eqref{def-D} for $D(s)$ to rewrite $\mathcal{E}_{2,k}$ as 
\begin{equation}\label{def2-bigSk} 
\mathcal{E}_{2,k}  = 
\sum_{{\begin{substack} {\varrho \in Z(\chi_{(k)}) \\ |\g| \ge k\g_0+t_0} \end{substack}}}^{\star}
\left( \e h(0) \mathcal{E}_{2,k,1}(\varrho) +\mathcal{E}_{2,k,2}(\varrho) \right),
\end{equation}
\begin{align*}
\text{ with }\ \mathcal{E}_{2,k,1}(\varrho) = &
 \frac{\sigma-\beta}{(\sigma-\beta)^2+\left(k \g_0 - \g\right)^2} + \frac{\sigma-1+\beta}{(\sigma-1+\beta)^2+\left(k \g_0 - \g\right)^2}
 \\& -  \k  \left( \frac{\sigma-\beta +\d }{(\sigma-\beta +\d )^2 + \left(k \g_0 - \g\right)^2}
+ \frac{\sigma-1+\beta +\d }{(\sigma-1+\beta +\d )^2 + \left(k \g_0 - \g\right)^2} \right),
\\
\mathcal{E}_{2,k,2}(\varrho)  = &   H(\sigma-\beta, k \g_0 - \g)  + H(\sigma-1+\beta,k \g_0 - \g)
\\&  - \k   \left( H(\sigma-\beta+\d , k \g_0 - \g) + H(\sigma-1+\beta+\d , k \g_0 - \g)  \right).
\end{align*}
The inequalities
\begin{align*}
 -\e \le &\sigma-\beta,\\
 1-2\e \le & \sigma-1+\beta +\d \le 1,\\
& \sigma-\beta+\d \le \d, \\
& \sigma-1+\beta +\d \le \d+1 ,
\end{align*}
and $|k\g_0-\g|  \ge t_0$ imply
\[
 \mathcal{E}_{2,k,1} (\varrho)
\ge \frac{1}{\left(k \g_0 - \g\right)^2} \left( -\e + \frac{1-2\e}{ t_0^{-2} +1} - \k    (2\d +1) \right).
\]
By Condition \eqref{choice-t0}, $\frac{1}{ 1+ t_0^{-2}} - \k    (2\d +1) >0$.
Thus
\begin{equation}\label{bound-bigSk1}
 \mathcal{E}_{2,k,1} (\varrho)
\ge - \frac{\e}{\left(k \g_0 - \g\right)^2} \left( 1 + \frac{2}{ t_0^{-2} +1} \right)
= - \frac{\e}{\left(k \g_0 - \g\right)^2} \frac{3+t_0^{-2}}{1+t_0^{-2}}.
\end{equation}
We use \eqref{boundF-H} and \eqref{boundF-m} 
to respectively bound $|H(\sigma-\beta,k \g_0 - \g )|$
and $|H(x, k \g_0 - \g)|$, with $x = \sigma-1+\beta,\sigma-\beta+\d $, or $\sigma-1+\beta+\d $.
Since  $-\frac{r}R \le  -\omega  \le \frac{\s-\b}\e$ and $x\ge \d-2\e_0$, 
we obtain
\begin{equation}\label{bound-bigSk2}
\left|\mathcal{E}_{2,k,2} (\varrho)\right|
\le \left(M(-r/R) \e^2 + \frac{1+2\k }{\d-2\e_0} m_h \e^3 \right)
\frac{1}{\left(\g - k\g_0\right)^2}.
\end{equation}
Together with \eqref{def2-bigSk}, \eqref{bound-bigSk1}, and \eqref{bound-bigSk2}, we deduce for $\mathcal{E}_{2,k} $:
\begin{equation}\label{bound-bigSk}
\mathcal{E}_{2,k} 
\ge \left[
 -  \left( \frac{3 + t_0^{-2}}{ 1+ t_0^{-2}} h(0) + M(-r/R)\right) 
-  \frac{1+2\k  }{\d-2\e_0}m_h  \e \right]
\e^2 \Sigma( k\g_0,t_0,\chi_{(k)}) ,
\end{equation}
\begin{equation}
\text{with }\ \Sigma(t,t_0,\chi) =
\sum_{{\begin{substack}{\varrho \in Z(\chi) \\ |\gamma| \ge t+t_0} \end{substack}}} \frac1{(\gamma-t)^2}.
\end{equation}
It remains to establish an upper bound for $\Sigma(t,t_0,\chi)$ when $t\ge 0$. In order to do this, we require explicit bounds for $N(T,\chi)$. 
Set
$\displaystyle{ 
\Phi(y) = \frac1{(y-t)^2} }$  and 
$\displaystyle{  \Phi'(y) = \frac{-2}{(y-t)^3}.
}$ 
By Stieltjes integration and Theorem \ref{bound-N}, we have
\[
 \Sigma(t,t_0,\chi) 
\le -\Phi(t+t_0)  N_2(t+t_0,q) -  \int_{t+t_0}^{\infty} \Phi'(y) N_1(y,q) dy.
\]
It remains to bound the last integral. 
\begin{lem}
\label{Lemme4} Let $t\ge0$ and $t_0$ be a positive integer.
Let $\chi$ be a primitive non-principal character modulo  $q$. Then
\[
\Sigma(0,t_0,\chi) 
\\ \le
(\log q) \left(\frac1{\pi t_0} + \frac{2 C_1}{ t_0^2}\right)  
+\frac{\log \frac{ e t_0 }{2\pi}}{\pi t_0} 
+  \frac{4 C_1  \log t_0 +C_1+4 C_2}{2 t_0^2} ,
\]
and when $t>0$,
\[
\Sigma(t,t_0,\chi) 
\le 
\frac1{\pi t_0} \log \frac{q(t+t_0)}{2\pi} + \frac{2 C_1}{ t_0^2}  \log (q(t+t_0))
+ \frac{\log\left( 1 + \frac{t}{t_0}\right)}{\pi t} 
\\ +  \frac{C_1}{ tt_0} \left( 1-\frac{t_0}{t} \log\left( 1 + \frac{t}{t_0}\right) \right)
+\frac{2 C_2}{ t_0^2}.
\]
\end{lem}
To bound $\e \Sigma(0,t_0,\chi_{k})$, we use $\frac{1}{r\log (q\g_0)} \le \frac{1}{r\log Q_0}$ and $\frac{\log q}{r\log (q\g_0)} \le \frac{\log q_1}{r\log Q_0}$, so that
\[
\e \Sigma(0,t_0,\chi) \le
\frac{ (\log q_1) \left(\frac1{\pi t_0} + \frac{2 C_1}{ t_0^2}\right)  
+\frac{\log \frac{ e t_0 }{2\pi}}{\pi t_0} 
+  \frac{4 C_1  \log t_0 +C_1+4 C_2}{2 t_0^2} }{r\log Q_0}.
\]
To bound $\e \Sigma(k\g_0,t_0,\chi_{k})$ when $k\ge1$, we also use
$\displaystyle{ 
\e \log\frac{q(k\g_0+t_0)}{2\pi} = \frac{\log \frac{q(k\g_0+t_0)}{2\pi}}{r\log (q\g_0)} 
\le \frac{\log  \frac{Q_0\left(k+\frac{t_0}{H_q}\right)}{2\pi} }{r\log Q_0}.
}$ 
In addition, we note that $\frac{\log\left( 1 + \frac{t}{t_0}\right)}{\pi t}  + \frac{C_1}{ tt_0} \left( 1-\frac{t_0}{t} \log\left( 1 + \frac{t}{t_0}\right) \right)$ is positive and decreases with $t$.
\begin{multline}
\text{Thus }\  \e \Sigma(k\g_0,t_0,\chi_{k})
\le 
\Big[\frac1{\pi t_0}  \log  \left ( \frac{Q_0}{2\pi}  \left(k+\frac{t_0}{H_q}\right)\right)
+ \frac{2 C_1}{ t_0^2} \log \left(Q_0\left(k+\frac{t_0}{H_q}\right) \right)
\Big.\\\Big. 
+ \frac1{\pi kH_q} \log\left( 1 + \frac{kH_q}{t_0}\right)
+  \frac{C_1}{ t_0 k H_q } \left( 1-\frac{t_0}{kH_q} \log\left( 1 + \frac{kH_q}{t_0}\right) \right)
+\frac{2 C_2}{ t_0^2}
\Big]\frac1{r\log Q_0}.
\end{multline}
We rearrange the terms and obtain $ \e \Sigma(k\g_0,t_0,\chi_{(k)}) \le s_k$, with
 \begin{align}
& \label{newdef-s0}
s_0 = \Big(   \frac{\log \Big(\frac{ e t_0 q_1 }{2\pi} \Big) }{\pi t_0} +  \frac{4 C_1  \log (q_1t_0)+C_1+4 C_2}{2 t_0^2} \Big) \e_0,
\\
&\label{newdef-sk}
s_k= \Big(
 \Big( \frac1{\pi t_0}  +\frac{2 C_1}{ t_0^2}\Big) \log  \Big( Q_0\Big(k+\frac{t_0}{H_q}\Big) \Big) 
- \frac{\log  (2\pi)}{\pi t_0}    
+\frac{2 C_2}{ t_0^2}
+ \frac{\Big(\frac1{\pi}-\frac{C_1}{kH_q}\Big) \log\Big( 1 + \frac{kH_q}{t_0}\Big)}{ kH_q} 
 +  \frac{C_1}{ t_0 k H_q }
\Big)  \e_0 ,
\end{align}
for $k\ge 1$
This allows to rewrite the bound \eqref{bound-bigSk} for $\mathcal{E}_{2,k}$ as
$\mathcal{E}_{2,k}  \ge- c_{2,2}(k,\e),$ 
with
\begin{equation}\label{newdef-c22}
  c_{2,2}(k,\e)
= \left(  \left( \frac{3 + t_0^{-2}}{ 1+ t_0^{-2}} h(0) + M(-r/R)\right) \e +  \frac{1+2\k  }{\d-2\e_0}m_h  \e^2  \right)s_k.
\end{equation}
This gives a final bound for $\mathcal{E}_2 $ and \eqref{ineq2-bigS} becomes
\begin{lem}\label{lem-mathcalE2}
Let $\s, \g_0, \k, \d, t_0$ be as in Section \ref{notations} and Proposition \ref{gal-Stechkin}.
Then
\begin{equation}\label{bound-mathcalE2}
- \mathcal{E}_2 
\le 
 - a_1 F(\sigma-\beta_0,0) 
 + \mathcal{C}_{2}(\e),
\end{equation}
\begin{equation}\label{newdef-C2}
\text{with }\ \mathcal{C}_{2}(\e) = a_1  c_{2,1}(\e) + \sum_{k=0}^{n_0}   a_k   c_{2,2}(k,\e),
\end{equation}
where $c_{2,1}$ and $c_{2,2}$ are defined in \eqref{newdef-c21} and \eqref{newdef-c22} respectively.
\end{lem}
\begin{rem}
Note that Lehman's method as used in \cite{MosTru} provides a larger bound with $\frac2{t_0^2}$ instead of $\frac{2C_1}{t_0^2}$ in the main factor $\log (q(t+t_0))$.
\end{rem}
\subsection{Study of the Gamma-term $\mathcal{E}_3$ 
}\label{Study-mathcalE3}
\begin{lem}\label{lem-mathcalE3}
Let $\s, \g_0, \k, \d, t_0$ be as in Section \ref{notations} and Proposition \ref{gal-Stechkin}.
Then
\begin{equation}\label{bound-mathcalE3}
\mathcal{E}_3
\le 
\frac{(1-\k) A h(0)}2 \e \log (q\g_0) + \mathcal{C}_{3}(\e),
\end{equation}
\begin{multline}\label{newdef-C3}
\text{ with }\  \mathcal{C}_{3}(\e) = \e h(0) \Bigg[ a_0 \left(-\frac{1-\k  }2\log \pi + \frac12 \frac{\Gamma'}{\Gamma} (3/2) - \frac{\k}2 \frac{\Gamma'}{\Gamma} ( (\s_0+2+\d)/2 )\right)
\Bigg. \\ \Bigg.  +  \sum_{k=1}^{n_0}  a_k \Big(
\frac{1-\k  }2 \log \frac{k}{2\pi}
+ \frac{\arctan \left(\frac{kH_q}{\s_0}\right) + \k   \arctan \left(\frac{kH_q}{\s_0+\d }\right)}{2kH_q} 
 + \frac{4+\k  (2+\d )^2}{4(kH_q)^2}
\Big)\Bigg].
\end{multline}
\end{lem}
\begin{proof}
$\displaystyle{\frac{\Gamma'}{\Gamma}}$ is an increasing function of the real variable, so for $k=0$ and $\s_0<\s<1$:
\begin{equation}\label{bound-Gamma1}
\mathfrak{G}_{\k,\d}(\s+2) =\frac{\Gamma'}{\Gamma}\left(\frac{\s+2}2\right) - \k \frac{\Gamma'}{\Gamma}\left(\frac{\s+2+\d }2\right)
\le \frac{\Gamma'}{\Gamma}\left(\frac{3}2\right) - \k \frac{\Gamma'}{\Gamma}\left(\frac{\s_0+2+\d }2\right). 
\end{equation}
For $k\ge 1$, 
the identity from \cite[page 12]{Mc2} gives
\begin{equation}\label{gamma-def-int}
\frac{\Gamma'}{\Gamma}\left(\frac{ x+iy}2\right)  = \frac12 \log
 \frac{x^2+y^2}{4} - \frac{x}{x^2+y^2} + \Re
\int_0^{+\infty}
\frac{u-[u]-\frac12 }{\left(u+\frac{x+iy}{2}\right)^2} d u.
\end{equation}
We have
$\displaystyle{
\Big| \Re \int_0^{+\infty}
\frac{u-[u]-\frac12 }{\big(u+\frac{x+iy}{2}\big)^2}d u \Big|
 \le \frac1y \arctan \left(\frac{y}{x}\right),
}$ 
so isolating $\log \frac{|y|}2$ in \eqref{gamma-def-int} gives the estimate
\begin{equation}\label{bound-gamma-def-int}
\left| \frac{\Gamma'}{\Gamma}\left(\frac{ x+iy}2\right)  - \left(\log \frac{|y|}2 -\frac{x}{x^2+y^2}\right) \right| \le \frac1y \arctan \left(\frac{y}{x}\right) 
+ \frac{x^2}{2y^2},
\end{equation}
\begin{multline}
\text{and then }\ 
\left|\mathfrak{G}_{\k,\d} ( \s+ik\g_0+\mathfrak{a} ) 
\right|
\le 
\left| \frac{1-\k}2 \log \frac{k\g_0}2 - \frac{\s+\mathfrak{a}}{(\s+\mathfrak{a})^2+(k\g_0)^2} + \k    \frac{\s+\mathfrak{a}+\d }{(\s+\mathfrak{a}+\d )^2+(k\g_0)^2}  \right|
\\ +  \frac{\arctan \left(\frac{k\g_0}{\s+\mathfrak{a}}\right) + \k \arctan \left(\frac{k\g_0}{\s+\mathfrak{a}+\d }\right) }{k\g_0} 
+ \frac{(\s+\mathfrak{a})^2 + \k (\s+\mathfrak{a}+\d )^2}{2(k\g_0)^2} .
\end{multline}
Both $\frac{x}{x^2+y^2}-\k   \frac{x+\d}{(x+\d)^2+y^2}$ and $\frac1y\arctan\frac{y}x$ are nonnegative and decreasing with $y$ when $\k < \frac{x}{x+\d}$ (which is the case from \eqref{def-kappa-delta} since $\k < \frac{\s_0}{1+\d} $).
Since $\g_0\ge H_q$ and $\s_0 < \s < 1$, then
\begin{equation}\label{bound-Gamma2}
\left|\mathfrak{G}_{\k,\d} ( \s+ik\g_0+\mathfrak{a} ) 
\right|
\le 
(1-\k  )\log \frac{k\g_0}2 
+ \frac{\arctan \left(\frac{kH_q}{\s_0}\right) + \k   \arctan \left(\frac{kH_q}{\s_0+\d }\right)}{kH_q} 
 + \frac{4+\k  (2+\d )^2}{2(kH_q)^2}.
\end{equation}
We conclude the proof by combining \eqref{def-mathcalE3} with \eqref{bound-Gamma1} and \eqref{bound-Gamma2}.
\end{proof}
\subsection{Study of the error term with the Gamma integral $\mathcal{E}_4$ 
}\label{Study-mathcalE4}
First, we establish a preliminary result to estimate $\Re \frac{\Gamma'}{\Gamma}\left(\frac{s}2\right) $:
\begin{lem}\label{lem-bound-Psi}
Let $\mathfrak{a}=0$ or $1$ and $T$ a real number. 
Let $u_{\mathfrak{a}}$ be the root of $\log\frac{u}2-\frac{\mathfrak{a}+1/2}{(\mathfrak{a}+1/2)^2+u^2}$,
\begin{align*}
& b_{\mathfrak{a}}  =\gamma + 3\log 2 + (-1)^{\mathfrak{a}}\frac{\pi}2  ,\\ 
& c_{\mathfrak{a}} =-\log2- \frac{\mathfrak{a}+1/2}{(\mathfrak{a}+1/2)^2+u_{\mathfrak{a}}^2} +  \frac1{u_{\mathfrak{a}}} \arctan \left(\frac{u_{\mathfrak{a}}}{\mathfrak{a}+1/2}\right) + \frac{\left(\mathfrak{a}+1/2\right)^2}{2u_{\mathfrak{a}}^2} .
\end{align*}
Then
\[ 
\left|\Re \frac{\Gamma'}{\Gamma}\left(\frac{ \mathfrak{a}+1/2+iT }2\right) \right|
\le U_{\mathfrak{a}}(T) =
\begin{cases}
b_{\mathfrak{a}}
&\textrm{if  }\ |T| \le u_{\mathfrak{a}},\\
\log |T|  +c_{\mathfrak{a}}
&\textrm{otherwise}.
\end{cases}\hfill 
\]
\end{lem}
Numerically, we have
\begin{table}[h]
\begin{tabular}{|l|l|l|l|}
\hline
 $\mathfrak{a}$ & $u_{\mathfrak{a}} $& $b_{\mathfrak{a}} $& $ c_{\mathfrak{a}} $\\
\hline
 $0$ & $2.2054\ldots$ & $4.2274 \ldots$ & $-0.1540\ldots$ \\
\hline
 $1$ & $2.4093\ldots$ & $1.0858\ldots$  & $-0.2647\ldots$ \\
\hline
\end{tabular}
\end{table} 
\begin{proof}
Thanks to Stirling's formula
$
\Re \frac{\Gamma'}{\Gamma}\left(\frac{  x+iy }2\right) = -\g -
\frac{2x}{x^2+y^2} + \sum_{n=1}^{+\infty} \left( \frac1n -
\frac{2(2n+x)}{(2n+x)^2+y^2} \right),
$
we see that $\Re \frac{\Gamma'}{\Gamma}\left(\frac{  \mathfrak{a}+1/2+iT  }2\right) $ increases with $|T|$. 
Thus when $|T|\le u_{\mathfrak{a}}$, 
\[
\left|\Re \frac{\Gamma'}{\Gamma}\left(\frac{  \mathfrak{a}+1/2+iT  }2\right) \right|
\le 
\left|\frac{\Gamma'}{\Gamma}\left(\frac{  \mathfrak{a}+1/2}2\right) \right|
= \gamma + 3\log 2 + (-1)^{\mathfrak{a}}\frac{\pi}2 .
\]
When $|T|\ge u_{\mathfrak{a}}$, we apply the inequality \eqref{bound-gamma-def-int} together with the fact that $\log\frac{|T|}2-\frac{\mathfrak{a}+1/2}{(\mathfrak{a}+1/2)^2+T^2}$ is nonnegative.
Bounding 
$
- \frac{\mathfrak{a}+1/2}{(\mathfrak{a}+1/2)^2+T^2} +  \frac1{|T|} \arctan \left(\frac{|T|}{\mathfrak{a}+1/2}\right) + \frac{\left(\mathfrak{a}+1/2\right)^2}{2T^2}
$
 with its value at $|T|=u_{\mathfrak{a}}$, it follows
\[
\left|\Re \frac{\Gamma'}{\Gamma}\left(\frac{  \mathfrak{a}+1/2+iT  }2\right) \right|
\le 
 \log \frac{|T|}2   - \frac{\mathfrak{a}+1/2}{(\mathfrak{a}+1/2)^2+u_{\mathfrak{a}}^2} 
 +  \frac{\arctan \left(\frac{u_{\mathfrak{a}}}{\mathfrak{a}+1/2}\right)}{u_{\mathfrak{a}}}  + \frac{\left(\mathfrak{a}+1/2\right)^2}{2u_{\mathfrak{a}}^2}.
\]
\end{proof}
Let $I_{\mathfrak{a}}$ be the integral
\begin{equation}\label{def-I}
I_{\mathfrak{a}}(x,y)=   \int_{0}^{+\infty}  \frac{ U_{\mathfrak{a}}(t)}{x(x^2+(y-t)^2)}  dt .
\end{equation}
Applying \eqref{boundF-m} to bound the $H$-terms in definition \eqref{def-mathcalE4}, it follows for $\s>\s_0$ that 
\begin{multline}
\label{bound1-mathcalE4}
\mathcal{E}_4
\le a_0  m_h \e^3\left(\frac1{\s_0^3} +  \frac{\k  }{(\s_0+\d)^3}
 + \frac{I_0(\s_0-1/2,0)+\k   I_0(\s_0+\d-1/2,0)}{\pi}  \right) 
\\
+ \frac{m_h \e^3}{2\pi} \sum_{k=1}^{n_0}  a_k
 \big[   I_{\mathfrak{a}}(\s_0-1/2,k\g_0) 
 + I_{\mathfrak{a}}(\s_0-1/2,-k\g_0) 
\big. \\  \big.
+ \k   \left( I_{\mathfrak{a}}(\s_0+\d-1/2,k\g_0)  + I_{\mathfrak{a}}(\s_0+\d-1/2,-k\g_0) \right) \big] .
\end{multline}
Note that $I_{\mathfrak{a}}(x,y)$ is decreasing with $x$.
The following lemma provides bounds for $I_{\mathfrak{a}}$.
\begin{lem}\label{lem-bound-Ia}
Let $\mathfrak{a}=0$ or $1$, $x>0,y>u_{\mathfrak{a}}+1$.
Then
\begin{align}
&\label{bound-I0a0}
 I_0(x,0) \le J_0(x),\\
& \label{bound-Ia-y}
 I_{\mathfrak{a}}(x,-y) \le J_{\mathfrak{a},1}(x,y) ,\\
& \label{bound-Iay}
I_{\mathfrak{a}}(x,y) \le J_{\mathfrak{a},2}(x,y) +J_{\mathfrak{a},3}(x,y) + (\log y) J_4(x,y), 
\end{align}
\begin{align*}
\text{with }\ 
& J_0(x)=  \frac{(b_0-c_0) \arctan\frac{u_0}{x}+   \frac{c_0\pi}{2}}{x^2} + \frac{\log u_0 +1}{x u_0}  ,
\\
& J_{\mathfrak{a},1}(x,y) 
= \frac{ (b_{\mathfrak{a}}-c_{\mathfrak{a}} ) \arctan\frac{y+u_{\mathfrak{a}}}{x}  - b_{\mathfrak{a}} \arctan\frac{y}{x} +  \frac{\pi c_{\mathfrak{a}}}{2}  }{x^2} 
-  \frac{\log (u_{\mathfrak{a}}+ y)}{xy} + \frac{u_{\mathfrak{a}}\log u_{\mathfrak{a}}}{xy(u_{\mathfrak{a}}+ y)} ,
\\
& J_{\mathfrak{a},2}(x,y) =  \frac{b_{\mathfrak{a}} \left( \arctan\frac{y}x  -\arctan\frac{y-u_{\mathfrak{a}}}x \right) }{x^2} ,
\\ 
& J_{\mathfrak{a},3}(x,y) = \frac{c_{\mathfrak{a}} \left( \frac{\pi}2 + \arctan\frac{y-u_{\mathfrak{a}}}{x} \right)}{x^2} ,  
\\ 
& J_4(x,y)= \frac{\pi }{x^2} + \frac{\pi }{x^2(\log y)^{3/2}}+ \frac{2\sqrt{\log y}}{xy}  + \frac{2}{xy(\log y)^{3/2}}.
\end{align*}
\end{lem}
\begin{proof}
Let $t\ge u_0$. We get \eqref{bound-I0a0} from the fact that $\frac{\log t }{x(x^2+t^2)} \le \frac{\log t }{xt^2} $, and thus
\[
I_0(x,0) 
 \le  \int_{0}^{u_0}   \frac{b_{0}}{x(x^2+t^2)}  dt + \int_{u_0}^{+\infty}   \frac{\log t }{xt^2}  dt + \int_{u_0}^{+\infty}   \frac{ c_0}{x(x^2+t^2)}  dt .
\] 
Let $t> u_{\mathfrak{a}}$. The announced inequality \eqref{bound-Ia-y} follows from $ \frac{\log t}{x(x^2+(y+t)^2)}\le  \frac{\log t}{x(y+t)^2}$ and from
\[
I_{\mathfrak{a}}(x,-y) 
 \le  \int_{0}^{u_{\mathfrak{a}}} \frac{b_{\mathfrak{a}}}{x(x^2+(y+t)^2)} dt 
+ \int_{u_{\mathfrak{a}}}^{+\infty} \frac{\log t}{x(y+t)^2} 
+ \frac{c_{\mathfrak{a}}}{x(x^2+(y+t)^2)} dt .
\]
Let $\ep>0$ be a parameter (depending on $y$), and split the integral $ I_{\mathfrak{a}}(x,y)$ at $t=u_{\mathfrak{a}}$.
The first integral is exactly $J_{\mathfrak{a},2}(x,y)$:
\begin{equation}\label{def-Ia2}
\int_{0}^{u_{\mathfrak{a}}} \frac{U_{\mathfrak{a}}(t)}{x(x^2+(y-t)^2)}  dt 
 = \frac{b_{\mathfrak{a}} }{x^2} \left( \arctan\frac{y}x  -\arctan\frac{y-u_{\mathfrak{a}}}x \right)
=J_{\mathfrak{a},2}(x,y). 
\end{equation}
Now for $t>u_{\mathfrak{a}} $, Lemma \ref{lem-bound-Psi} gives $U_{\mathfrak{a}}(t) = \log t + c_{\mathfrak{a}}$.
It is immediate that
\begin{equation}\label{def-Ia3}
\int_{u_{\mathfrak{a}}}^{+\infty} \frac{ c_{\mathfrak{a}}}{x(x^2+(y-t)^2)} dt 
=\frac{c_{\mathfrak{a}}}{x^2} \left( \frac{\pi}2 + \arctan\frac{y-u_{\mathfrak{a}}}{x} \right)
=J_{\mathfrak{a},3}(x,y) .
\end{equation}
To compute the $(\log t)$-part, we split at $t=y(1-\epsilon)$ and $y(1+\epsilon)$.\\
When $u_{\mathfrak{a}} < t< y(1-\epsilon)$ or $t > y(1-\epsilon)$, we use $ \tfrac{\log t}{x^2+(y-t)^2}\le  \tfrac{\log t}{(y-t)^2}$. \\
Assuming $y>u_{\mathfrak{a}}+1$, we obtain
\begin{align*}
& \left(\int_{u_{\mathfrak{a}}}^{y(1-\epsilon)} +\int_{y(1+\epsilon)}^{+\infty} \right) \frac{\log t}{x(y-t)^2} dt 
\\ &= 
\frac1{xy\ep}
\Big(
 -\ep \log (y-u_{\mathfrak{a}}) + \frac{ \ep u_{\mathfrak{a}}\log u_{\mathfrak{a}}}{u_{\mathfrak{a}} - y} 
+ 2\log y 
 + (1-\epsilon)\log (1-\epsilon)
+ (1+\epsilon)\log (1+\epsilon)
\Big)
\\ &\le 
\frac{2\left( \log y  +\epsilon^2 \right)}{xy\ep}.
\end{align*}
When $y(1-\epsilon) < t< y(1+\epsilon)$, we have 
\[
\int_{y(1-\epsilon)}^{y(1+\epsilon)}  \frac{\log t}{x(x^2+(y-t)^2)} dt 
 \le  \int_{y(1-\epsilon)}^{y(1+\epsilon)}  \frac{\log (y(1+\epsilon))}{x(x^2+(y-t)^2)} dt 
 =  \frac{2 \log (y(1+\epsilon)) \arctan\frac{y\epsilon}{x}  }{x^2} 
 \le \frac{\pi(\log y + \ep)}{x^2} .
\]
Choosing $\ep=\frac{1}{ \sqrt{ \log y } }$ and adding the two last inequalities, we obtain
\begin{equation}\label{def-Ia4}
\int_{u_{\mathfrak{a}}}^{+\infty}  \frac{\log t}{x(x^2+(y-t)^2)} dt 
\le 
(\log y)   J_4(x,y).
\end{equation}
\eqref{bound-Iay} follows from \eqref{def-Ia2}, \eqref{def-Ia3} and \eqref{def-Ia4}.
\end{proof}
We combine the bound \eqref{bound1-mathcalE4} for $\mathcal{E}_4$ with Lemma \ref{lem-bound-Ia}.
In particular, to bound $\e(I_{\mathfrak{a}}(\s_0-1/2,k\g_0) + \k I_{\mathfrak{a}}(\s_0+\d-1/2,k\g_0))$, we use the fact that 
$
\e \log(k\g_0) 
= \frac{1+\frac{\log k}{\log \g_0}}{r(1+\frac{\log q}{\log \g_0})}
 \le 
\frac{\log (kH_q)}{r\log H_q} .
$
\begin{lem}\label{lem-mathcalE4}
Let $\s, \g_0, \k, \d, t_0$ be as in Section \ref{notations} and Proposition \ref{gal-Stechkin}.
Then
\begin{equation}\label{bound-mathcalE4} 
 \mathcal{E}_4 
\le \mathcal{C}_{4}(\e) ,
\end{equation}
\begin{multline}\label{newdef-C4}
\text{with }\ 
\mathcal{C}_{4}(\e) = 
 \e^3 a_0 m_h \left( \frac1{\s_0^3}+\k   \frac1{(\s_0+\d)^3} + \frac{J_0(\s_0-1/2) + \k   J_0(\s_0-1/2+\d)}{\pi} \right) 
\\ + \e^3 \frac{m_h }{2\pi}  \sum_{k=1}^{n_0} \sum_{j=1}^3   a_k  \max_{\mathfrak{a}=0,1} 
 \left( J_{\mathfrak{a},j}(\s_0-1/2,kH_q)  + \k   J_{\mathfrak{a},j}(\s_0-1/2+\d,kH_q) \right)
\\
+ 
\e^2 \frac{m_h }{2r\pi} \sum_{k=1}^{n_0}   a_k  \left( J_4(\s_0-1/2,kH_q) + \k   J_4(\s_0-1/2+\d,kH_q) \right).
\end{multline}
\end{lem}
\subsection{Choice of the smooth weight $h$}
\label{choice-weight}
Let $\l>0$ and $c_1,c_2,x_1,x_2$ some fixed real numbers. 
Consider
\begin{equation}\label{def-g5}
g(x) = c_1 \cos(x_1 x)+c_2 \cosh(x_2 x) - (c_1\cos(x_1 \l)+c_2\cosh(x_2 \l))
\end{equation}
for $-\l \le x \le \l$ and $g(x) = 0$ otherwise.
We have
$
g^{\prime\prime}(x) = -c_1\cos(x_1u) x_1^2 + c_2 \cosh(x_2 u) x_2^2.
$
Note that $g$ is compactly supported, non-negative, even, infinitely differentiable, and that $g^{\prime\prime}$ is even, negative, and increasing on $(0,\l)$.
Let 
$ d_5=2\l,$
and define $h$ the self-convolution of $g$:
\[
h(u) = (g\star g) (u) = \int_{u-\l}^{\l} g(x) g(u-x) d x, 
\]
when $0< u < d_5$, and $h(u) = 0$ otherwise.
Note that $h$ satisfies Conditions \eqref{cond-f}.
We have
$\displaystyle{
h^{\prime\prime}(u) = \int_{-\infty}^{\infty} g(x) g^{\prime\prime}(u-x) dx
}$
and denote 
$
m_h = \max_{0< u < d_5}\left| h^{\prime\prime}(u)\right|.
$
Numerically, we observe that $m_h = - h^{\prime\prime}(0)$. 
\newline
\begin{rem}
In his famous article about the least prime in an arithmetic progression \cite[Lemmas 7.1-7.4]{HB}, Heath-Brown used various smoothed versions of $-\frac{L'}{L}(s,\chi)$. 
Later Xylouris followed one of Heath-Brown's remarks and proposed some other weights. 
For instance \cite{Kad1} used the weight as defined in \cite[Lemma 7.4]{HB}.
In \cite{JanKwo} Jan and Kwon compared all the weights proposed by Heath-Brown and Xylouris which allowed them to improve the last zero-free region for zeta (by $1.322\%$). 
After investigating all these five weights, we find that the weight given in \cite[page 75]{Xyl} by Xylouris provides here the best constant for the zero-free region.
\end{rem}
\subsection{Choice of the trigonometric polynomial $P(x)$}
\label{choice-pol-trig}
The last area where we seek some final improvement concerns the trigonometric polynomial. 
We take here the opportunity to briefly describe what we know about the progresses concerning this aspect.
For the Riemann zeta function, de la Vall\'ee Poussin used the following quadratic polynomial:
$\displaystyle{
2(1+\cos x)^2 = 3+4\cos x + \cos(2x),
}$
and Landau showed that the value for the constant $R_0$ was given by
\begin{equation}
\label{def-R0-theory}
\frac{A}{2\left(\sqrt{a_1}-\sqrt{a_0}\right)^2}, 
\end{equation}
with $A=\sum_{k=1}^{n_0}  a_k$.
In the case of de la Vall\'ee Poussin's quadratic polynomial, \eqref{def-R0-theory} equals $ 34.82050\ldots$.
Finding a polynomial satisfying \eqref{cond-pol-trig} and such that the value for \eqref{def-R0-theory} is as small as possible becomes challenging as the degree of the polynomial increases.
In analysis, this problem is referred to as Landau's extremal problems (see \cite[Section 9]{Rev}). 
In 1970, Stechkin \cite{St} introduced the degree 4 polynomial:
\begin{multline}\label{def-Stechkin-trig-pol}
P(x)
= 8(0.9126+\cos x)^2(0.2766+\cos x)^2 
\\ = 
11.18...
+ 19.07... \cos x
+ 11.67... \cos (2x)
+ 4.75... \cos (3x)
+ \cos (4x),
\end{multline}
which brings \eqref{def-R0-theory} down to $17.42622\ldots$.
Together with his clever idea mentioned in Section \ref{section-kappa}, he was able  to reduce the constant in the zero-free region to $9.65$.
In 1975, minor modifications of Stechkin's idea led Rosser and Schoenfeld \cite{RS2} to compute $9.645908801$ instead, and in 1984, McCurley \cite{Mc2} generalized this to Dirichlet $L$-functions. 
In earlier work \cite{Kad1} and in \cite{Kad0}, we used a minor modification of Stechkin's polynomial \eqref{def-Stechkin-trig-pol}:
\begin{multline}\label{def-Kadiri-trig-pol}
P(x)
= 8(0.91+\cos x)^2(0.265+\cos x)^2 
\\ = 
10.91...
+ 18.63... \cos x
+ 11.45... \cos (2x)
+ 4.7 \cos (3x)
+ \cos (4x).
\end{multline}
We reduced the constant for zeta to $5.69693$ \cite{Kad1}, which was the first significant result 35 years after Stechkin, 
as well as the constant for Dirichlet $L$-functions of McCurley's \cite{Mc2}, obtaining the constant $6.4355$ \cite{Kad0}.
In 2012, I first heard through Olivier Ramar\'e about a forgotten article of Kondratev: 
in \cite{Kon} he found a polynomial of degree 8 which lead to a constant $17.27230\ldots$ for \eqref{def-R0-theory}.
Consequently, it reduced Stechkin's constant for zeta's zero-free region to $9.54789695$.
In 2014, Mossinghoff and Trudgian \cite{MosTru} used the author's method and also investigated numerically other trigonometric polynomials of higher degrees (up to 40).
Among these, they find the following polynomial of degree 16 to produce the smallest constant for zeta's zero-free region: 
\begin{small}
\begin{table}[ht] 
\caption{$P(x) = \left|\sum_{k=0}^{16} c_k e^{ikx} \right|^2 =   \sum_{k=0}^{16} a_k \cos(kx)$}
\label{trig-pol-16}
\begin{tabular}{| l l | l l |}
\hline
$c_0$ & 1 & $a_0$ & 1 \\
$c_1$ & $-2.09100370089199$ & $a_1$ & $1.74126664022806$ \\
$c_2$ & $0.414661861733616$ & $a_2$ &  $1.128282822804652 $ \\
$c_3$ & $4.94973437766435$ & $a_3$ &  $0.5065272432186642 $ \\
$c_4$ & $2.26052224951171$ & $a_4$ &  $0.1253566902628852 $ \\
$c_5$ &  $8.58599241204357 $ & $a_5$ &  $9.35696526707405\cdot10^{-13} $ \\
$c_6$ &  $6.87053689828658 $ & $a_6$ &  $4.546614790384321\cdot10^{-13} $ \\
$c_7$ &  $22.6412990090005 $ & $a_7$ &  $0.01201214561729989 $ \\
$c_8$ &  $6.76222005424994 $ & $a_8$ &  $0.006875849760911001 $ \\
$c_9$ &  $50.2233943767588 $ & $a_9$ &  $7.77030543093611\cdot10^{-12} $ \\
$c_{10}$ &  $8.07550113395201 $ & $a_{10}$ &  $2.846662294985367\cdot10^{-7} $ \\
$c_{11}$ &  $223.771572768515 $ & $a_{11}$ &  $0.001608306592372963 $ \\
$c_{12}$ &  $487.278135806977 $ & $a_{12}$ &  $0.001017994683287104 $ \\
$c_{13}$ &  $597.268928658734 $ & $a_{13}$ &  $2.838909054508971\cdot10^{-7} $ \\
$c_{14}$ &  $473.937203439807 $ & $a_{14}$ &  $ 5.482482041999887\cdot10^{-6} $ \\
$c_{15}$ &  $237.271715181426 $ & $a_{15}$ &  $ 2.412958794855076\cdot10^{-4}$ \\
$c_{16}$ &  $59.6961898512813 $ & $a_{16}$ &  $ 1.281001290654868\cdot10^{-4}$ \\
 &   & $A$ &  $3.523323140225021$ \\
\hline 
\end{tabular}
\end{table}
\end{small}\ \\
We note that this polynomial is close to best possible:
with it \eqref{def-R0-theory} equals $17.24998\ldots$ while the optimal value for \eqref{def-R0-theory} over all even trigonometric polynomial with non-negative coefficients and satisfying $a_1>a_0$ is no smaller than $17.23415$ (as proven in \cite[Theorem 2]{MosTru}).
\subsection{Computations}
\label{Computations}
We use the polynomial given in Section \ref{choice-pol-trig} and the smooth function $h$ depending on $c_1,c_2,x_1,x_2,\lambda$ as defined in Section \ref{choice-weight}.
We set $c_1=1$ and $x_1 = 1$ and for each set values for $c_2,x_2,$ and $\lambda$, we choose $r$ and $\eta_1$ so as to make the zero-free region constant $R_0$ as given by \eqref{formula-R0} as small as possible. 
We repeat the step with replacing the value of $R$ at step $k$ by the value of $R_0$ at step $k+1$.
We repeat until we get the first two decimal digits between $R$ and $r$ to be the same. 
We record our results in the following tables and display up to the first four decimal digits for each value computed (that is for $\kappa, \delta$, and $R_0$).
\begin{small}
\begin{table}[ht] 
\caption{For $3\le q \le 1\,000$}\label{table-R01}
\begin{tabular}{|c|c|c|c|c|c|c|c|c|c|}
\hline Step & $R$         & $r$      & $c_2$   & $x_2$   & $\lambda$ & $\eta_1$ & $\kappa$     & $\delta$      & $R_0$ \\
\hline
\hline 1 & $9.6459$ & $5.857$ & $1.361$ & $0.765$ & $0.551$ & $0.0007$ & $0.4414$ & $0.6198$ & $5.8579$\\
\hline 2 & $5.8579$ & $5.622$ & $1.875$ & $0.639$ & $0.533$ & $0.0011$ & $0.4381$ & $0.6208$ & $5.6223$\\
\hline 3 & $5.6223$ & $5.599$ & $0.194$ & $1.650$ & $0.531$ & $0.0012$ & $0.4378$ & $0.6209$ & $5.5992$\\
\hline 4 & $5.5992$ & $5.596$ & $0.111$ & $2.013$ & $0.531$ & $0.0012$ & $0.4377$ & $0.6210$ & $5.5968$\\
\hline
\end{tabular}
\end{table}
\begin{table}[ht] 
\caption{For $1000< q \le 400\,000$}\label{table-R02}
\begin{tabular}{|c|c|c|c|c|c|c|c|c|c|}
\hline Step & $R$         & $r$      & $c_2$   & $x_2$   & $\lambda$ & $\eta_1$ & $\kappa$     & $\delta$      & $R_0$ \\
\hline
\hline 1 & $9.6459$ & $5.857$ & $1.384$ & $0.759$ & $0.551$ & $0.0007$ & $0.4414$ & $0.6198$ & $5.8579$ \\
\hline 2 & $5.8579$ & $5.622$ & $1.494$ & $0.710$ & $0.533$ & $0.0011$ & $0.4381$ & $0.6208$ & $5.6223$ \\
\hline 3 & $5.6223$ & $5.599$ & $0.189$ & $1.665$ & $0.531$ & $0.0012$ & $0.4378$ & $0.6209$  & $5.5992$\\
\hline 4 & $5.5992$ & $5.596$ & $0.086$ & $2.190$ & $0.531$ & $0.0012$ & $0.4377$ & $0.6210$  & $5.5968$\\
\hline
\end{tabular}
\end{table}
\end{small}
\section{Some complementary proofs}\label{Complementary-Proofs}
%
\subsection{Explicit formulae for Dirichlet $L$-functions}
We give here an explicit formula relating sums over zeros of a Dirichlet $L$-function and sums over primes.
This corrects a mistake in \cite[Theorem 3.1]{Kad1} for the case of primitive characters.
\begin{thm}
\label{FormuleExplicite}
Let $\phi$ be a complex valued function such that 
\begin{enumerate}
\item[(A)] $\phi$ is $C^1$ on $\mathbb{R}-S$, where $S$ is a finite set of points $a_i$ where both $\phi$ and its derivative $\phi'$ have at worse removable discontinuities.
Moreover, at these points $\phi$ verifies $\phi(a_i) = \frac12 [\phi(a_i+0)+\phi(a_i-0)]$.
\item[(B)] There exists $b>0$ such that $\phi(x)e^{x/2}$ and $\phi'(x)e^{x/2}$ are
$\mathcal{O}(e^{-(\frac12 +b)|x|})$ as $x\to\infty$.
\end{enumerate}\noindent
For each $a<1$ verifying $0<a<b$, $\phi$ has a Laplace transform
$
\Phi (s) = \int_{0}^{+\infty}\phi (x)e^{-sx}d x
$
which is holomorphic in $-(1+a)<\sigma<a$ and which is $\mathcal{O}(1/|t|)$ in $-(1+a) \le \sigma \le a$.
Let $q\in \N$ and $\chi$ a primitive character modulo $q$.
Let $\mathfrak{a}= 0$ if $\chi(-1)=1$, $1$ if $\chi(-1)=-1$.
Then
\begin{multline}\label{explicit-formula-phi-primitive}
 \sum_{n\ge1}\Lambda(n)\chi(n)\phi(\log n)
=
 \phi(0)\log \frac q\pi 
- \sum_{\varrho \in Z(\chi)}  \Phi(-\varrho)
- \sum_{n\ge1}\frac{\Lambda(n)\overline\chi(n) \phi(-\log n)}{n}
\\ +\frac1{2i\pi} \int_{\frac12-i\infty}^{\frac12+i\infty}\Re \frac{\Gamma'}{\Gamma} \left(\frac{s+\mathfrak{a}}2\right) \Phi(-s) d s,
\end{multline}
\begin{multline}\label{explicit-formula-phi-zeta}
\text{and}\ 
 \sum_{n\ge1}\Lambda(n) \phi(\log n)
=
\Phi(-1)+ \Phi(0) -\phi(0)\log \pi  - \sum_{\varrho \in Z(\zeta)}  \Phi(-\varrho) 
 - \sum_{n\ge1}\frac{\Lambda(n) \phi(-\log n)}{n}
\\+ \frac1{2i\pi} \int_{\frac12-i\infty}^{\frac12+i\infty}\Re \frac{\Gamma'}{\Gamma} \left(\frac{s}2\right) \Phi(-s) d s,
\end{multline}
where $Z(\chi)$ and $Z(\zeta)$ are the sets of non-trivial zeros of respectively $\L(s,\chi)$ and $\zeta(s)$.
\end{thm}
\begin{proof}
The Inverse Laplace transform gives
$
\phi(\log n) = \frac1{2i\pi} \int_{-(1+a)-i\infty}^{-(1+a)+i\infty}\Phi(s) n^s d s \quad(n\ge1).
$
Combining this with a change of variable ($s\to -s$)
, we obtain
\begin{equation}\label{eval-Dirichlet-sum}
 \sum_{n\ge1}\Lambda(n)\chi(n)\phi(\log n)
=\frac1{2i\pi} \int_{1+a-i\infty}^{1+a+i\infty} -\frac{L'}{L}(s,\chi) \Phi(-s) d s
= \mathcal{I}_{1+a-i\infty}^{1+a+i\infty},
\end{equation}
with
$\mathcal{I}_{x+iy}^{x'+iy'} = \frac1{2i\pi} \int_{x+iy}^{x'+iy'} -\frac{L'}{L}(s,\chi) \Phi(-s) d s$.
We consider $\mathcal{I}_{1+a-iT}^{1+a+iT}$ the truncated integral, where $T>0$ is fixed and does not equal the ordinate of a zero of $L(s,\chi)$.
We move the contour of integration to $[-a-iT,-a+iT]$, and collect the residues of $-\frac{L'}{L}(s,\chi) \Phi(-s)$ in the issued rectangle.
The poles in this rectangle 
are simple and located at the non-trivial zeros $\varrho\in Z(\chi)$, each with residue $-1$.
In addition, if $\chi$ is even, there is another simple pole at $s=0$, also with residue $-1$.
Letting $\mathcal{P}_{T,\chi}$ be the pole contribution in the rectangle, we have from Cauchy's Residue Theorem that 
\begin{equation}\label{Cauchy}
\mathcal{I}_{1+a-iT}^{1+a+iT}
= \mathcal{P}_{T,\chi} + \mathcal{I}_{-a-iT}^{-a+iT} + \mathcal{I}_{-a+iT}^{1+a+iT} - \mathcal{I}_{-a-iT}^{1+a-iT},
\end{equation}
with
\begin{equation}\label{Poles-primitive}
 \mathcal{P}_{T,\chi} = - (1-\mathfrak{a}) \Phi(0) - \sum_{\begin{substack} {\varrho \in Z(\chi) \\ |\Im \varrho|<T} \end{substack}}  \Phi(-\varrho). 
\end{equation}
It follows from Condition (B) that $\mathcal{I}_{-a+iT}^{1+a+iT}$ and $\mathcal{I}_{-a-iT}^{1+a-iT}$ have limit $0$ as $T\to \infty$.
Moreover, the functional equation \cite[(13) (14) page 71]{Dav}
\begin{equation}\label{functional-equation-Dirichlet}
 -\frac{L'}{L}(s,\chi) = \log \frac q\pi + \frac{L'}{L}(1-s,\overline\chi) + \frac12
\Big\{ \frac{\Gamma'}{\Gamma}\left(\frac{s+\mathfrak{a}}{2}\right)
+\frac{\Gamma'}{\Gamma}\left(\frac{1-s+\mathfrak{a}}{2}\right) \Big\}
\end{equation}
allows to split $ \mathcal{I}_{-a-i\infty}^{-a+i\infty}$ as the sum of three integrals:
\begin{equation}\label{eval-I}
 \mathcal{I}_{-a-i\infty}^{-a+i\infty} =  I_{ \mathfrak{a},1}(\chi) + I_{ \mathfrak{a},2}(\chi) + I_{ \mathfrak{a},3}(\chi) , 
\end{equation}
\begin{align}
\text{with}\ & \label{eval-I1}
I_{ \mathfrak{a},1}(\chi) = \frac1{2i\pi} \int_{-a-i\infty}^{-a+i\infty}  \log \frac q\pi  \Phi(-s) d s=  \phi(0)\log \frac q\pi , 
\\
& \label{eval-I2}
I_{ \mathfrak{a},2}(\chi) = \frac1{2i\pi} \int_{-a-i\infty}^{-a+i\infty} \frac{L'}{L}(1-s,\overline\chi) \Phi(-s) d s = -\sum_{n\ge1}\frac{\Lambda(n)\overline\chi(n) \phi(-\log n)}{n}, 
\\
&
I_{ \mathfrak{a},3}(\chi) =  \frac1{2i\pi} \int_{-a-i\infty}^{-a+i\infty} \frac12 \left( \frac{\Gamma'}{\Gamma}\left(\frac{s+\mathfrak{a}}{2}\right) +\frac{\Gamma'}{\Gamma}\left(\frac{1-s+\mathfrak{a}}{2}\right) \right) \Phi(-s) d s.
\end{align}
To evaluate $I_{ \mathfrak{a},3}(\chi)$, we move the path of integration to the $\frac12 $-line on which $\Gamma$ verifies 
\[
\frac12 \left( \frac{\Gamma'}{\Gamma}\left(\frac{s+\mathfrak{a}}{2}\right) +\frac{\Gamma'}{\Gamma}\left(\frac{1-s+\mathfrak{a}}{2}\right) \right)
=\Re \frac{\Gamma'}{\Gamma}\left(\frac{s+\mathfrak{a}}{2}\right).
\]
Condition (B) and the fact that $\frac{\Gamma'}{\Gamma}(s)$ has a simple pole at $s=0$ with residue $-1$
lead to
\begin{equation}  \label{eval-I3}
I_{ \mathfrak{a},3}(\chi)
=\frac1{2i\pi} \int_{\frac12-i\infty}^{\frac12+i\infty}\Re \frac{\Gamma'}{\Gamma} ( s+\mathfrak{a} ) \Phi(-s) d s 
+ (1-\mathfrak{a}) \Phi(0).
\end{equation}
Note that this contribution of $\Phi(0)$ cancels the one arising in \eqref{Poles-primitive}.
Finally, together with \eqref{Cauchy}, \eqref{eval-I}, \eqref{eval-I1}, \eqref{eval-I2}, and \eqref{eval-I3}, we rewrite \eqref{eval-Dirichlet-sum} in the form announced in \eqref{explicit-formula-phi-primitive}.
\\
The formula \eqref{explicit-formula-phi-zeta} is obtained similarly. 
We recall that the poles of $-\frac{\zeta'}{\zeta}(s)$ are all simple and are located at $s=1$, with residue $1$, and at the non-trivial zeros $\varrho\in Z(\zeta)$, with residue $-1$.
In this case the polar contribution is
$
 \mathcal{P}_{T} = \Phi(-1) - \sum_{\begin{substack} {\varrho \in Z(\zeta) \\ |\Im \varrho|<T} \end{substack}}    \Phi(-\varrho),
$
instead of \eqref{Poles-primitive}, and 
we use the functional equation 
$
 -\frac{\zeta'}{\zeta}(s) = -\log \pi + \frac{\zeta'}{\zeta}(1-s) 
+ \frac12 \left( \frac{\Gamma'}{\Gamma}\left(\frac{s}{2}\right) +\frac{\Gamma'}{\Gamma}\left(\frac{1-s}{2}\right) \right)
$
of \cite[page 59]{Dav} 
instead of \eqref{functional-equation-Dirichlet}.
\end{proof}
\subsection{Proof of Proposition \ref{Prop1.1}}\label{proof-ExplicitFormula}
Let $s\in \C$.
We consider $ \phi(y) = (f(0)-f(y))e^{-ys} $ if $y\ge0$, and $ \phi(y) =0$ otherwise.
Thus $\phi$ verifies all the conditions stated in Theorem \ref{FormuleExplicite}.
We have $\phi(0)=0$, $\phi(-\log n)=0$, and for each $\Re z < \Re s$, 
$
\Phi(-z) = \frac{f(0)}{s-z} - F(s-z) =-\frac{F_2(s-z)}{(s-z)^2}.
$
We insert this definition in the explicit formula \eqref{explicit-formula-phi-primitive}, and take its real part.
Together with the classical explicit formula \cite[Equations (17) (18) on p 83]{Dav} for $-\Re \frac{L'}{L}(s,\chi)$, namely
\begin{equation}\label{classic-formula-L}
\Re \sum_{n\ge1}\frac{\Lambda(n)\chi(n) }{n^s} 
= \frac12 \log \frac{q}{\pi} + \frac12\Re \frac{\Gamma'}{\Gamma} \left(\frac{ s+\mathfrak{a}}2\right)  - \sum_{\varrho \in Z(\chi)} \Re \frac1{s-\varrho},
\end{equation}
\begin{multline*} 
\text{we obtain }\ f(0)\Big(\frac12 \log \frac{q}{\pi} + \frac12\Re \frac{\Gamma'}{\Gamma} \left(\frac{ s+\mathfrak{a}}2\right) - \sum_{\varrho \in Z(\chi)} \Re \frac1{s-\varrho} \Big)
- \Re \sum_{n\ge1}\frac{\Lambda(n)\chi(n) }{n^s}f(\log n) 
\\= 
- f(0) \sum_{\varrho \in Z(\chi)}  \Re \frac1{s-\varrho}
+ \sum_{\varrho \in Z(\chi)} \Re F(s-\varrho)
- \Re \frac1{2i\pi} \int_{\frac12-i\infty}^{\frac12+i\infty}\Re \frac{\Gamma'}{\Gamma}\left(\frac{ z+\mathfrak{a}}2\right)\frac{F_2(s-z)}{(s-z)^2} d z.
\end{multline*}
This establishes \eqref{Explicit-formula-primitive-f} for $\Re s>1$.
Since the functions defined on the left and right hand side of the equality are both  harmonic functions defined on the whole complex plane, then the identity extends to all $s\in \C$.

We prove \eqref{Explicit-formula-zeta-f} in a similar manner: we use \eqref{explicit-formula-phi-zeta} with 
$ \Phi(-1)=  \frac{f(0)}{s-1} - F(s-1), \Phi(0)= -\frac{F_2(s)}{s^2},$ 
and the classical explicit formula 
\cite[(8) (11) page 80]{Dav} 
\begin{equation}\label{classic-formula-zeta}
\Re \sum_{n\ge1}\frac{\Lambda(n)}{n^s}
=  \frac1{s-1} - \frac12 \log \pi + \frac12\Re \frac{\Gamma'}{\Gamma} \left(\frac{s+2}{2}\right)  - \sum_{\varrho \in Z(\zeta)} \Re \frac1{s-\varrho}.
\end{equation}
\subsection{Handling the non-primitive characters} 
We introduce similar notation to \cite[Section 3.4.]{RR}:
\begin{equation}\label{def-cp}
c_p(\s,\k,\d) = \sum_{m\ge1} \frac{1}{p^{m\s}} \left(1-\frac{\k}{p^{m\d}}\right) 
= \frac1{p^{\s}-1}-\frac{\k}{p^{\s+\d}-1}. 
\end{equation}
\begin{lem}\label{primitivity}
Let $\s, \g_0, \k, \d, t_0$ be as in Section \ref{notations} and Proposition \ref{gal-Stechkin}.
Let  $n_0 \ge 2$ and $P(x) = \sum_{k=0}^{n_0}a_k$ be the trigonometric polynomial as in Section \ref{choice-pol-trig}.
In addition we assume that these coefficients and $\s_0,\k,\d$ satisfy \eqref{cond-primitivity-kappa-delta2}
\begin{equation}
a_0 +  \left( \frac{1-\kappa}{\max\left(c_2(\s_0,\k,\d),2c_3(\s_0,\k,\d)\right)} - 1 \right) \sum_{k=2}^{n_0}  a_k  >0 .
\end{equation}
Then
\begin{equation}\label{thesum1}
 \frac{1-\kappa}2 f(0) \sum_{k=1}^{n_0}  a_k \log \left(\frac{q}{q_k}\right)
+ \sum_{k=0}^{n_0}  a_k S(\sigma+ik\gamma_0,\chi_{(k)}-\chi^k)
\ge 0.
\end{equation}
\end{lem}
\begin{proof}
We input the definition of $S(s,\chi)$ in \eqref{thesum1}. Note that for $k=1$ the factor of $a_1$ vanishes.
The left term equals
\begin{multline}\label{inq}
a_0 \sum_{p|q} \sum_{m\ge 1}  \frac{\log p}{p^{m\s}}f(m\log p) \left(1-\frac{\k}{p^{m\d}}\right) 
+ \sum_{k=2}^{n_0}  a_k \Bigg( \frac{1-\kappa}2 f(0) \log \frac{q}{q_k}
\\ + \sum_{\begin{substack}{p|q \\p \nmid q_k}\end{substack}} \sum_{m\ge 1}  \Re\left(\frac{\chi_{(k)}(p^m)}{p^{imk\g_0}}\right)\frac{(\log p)f(m\log p)}{p^{m\s}} \left(1-\frac{\k}{p^{m\d}}\right) \Bigg) .
\end{multline}
For $k\ge 2$, we use the inequalities
$\log \frac{q}{q_k} \ge \sum_{\begin{substack}{p|q\\ p \nmid q_k}\end{substack}} \nu_p\left(\frac{q}{q_k}\right)\log p,$
$f(0) \ge \sum_{m\ge1} \frac{f(m\log p)}{p^{m\s}c_p(\s,\k,\d)  } \left(1-\frac{\k}{p^{m\d}}\right)$,
and
$\Re \left(\frac{\chi_{(k)}(p^m)}{p^{imk\gamma_0}}\right) \ge-1$,
where $\nu_p$ is the notation for $p$-adic valuation.
Using the notation
\[
C_p(\s,\k,\d)  = a_0 +  \sum_{k=2}^{n_0}  a_k \left( \frac{1-\kappa}{2c_p(\s,\k,\d)}\nu_p\left(\frac{q}{q_k}\right) - 1 \right) >0 ,
\]
we have that \eqref{inq} is larger than 
\[
\sum_{p|q} \sum_{m\ge 1} \frac{(\log p) f(m\log p)}{p^{m\s}} \left( 1-\frac{\k}{p^{m\d}} \right) C_p(\s,\k,\d).
\]
It is immediate that $c_{p}(\s)$ is positive and increases as $\d$ increases and decreases as $\k$ increases. It also it decreases as $p$ or $\s$ increases: 
this is easily verified by showing that the respective derivatives have opposite sign of
$\s \left(1-p^{-(\s+\d)}\right)^2 - \k (\s+\d) p^{-\d} \left(1-p^{-\s}\right)^2$, 
which is positive under the conditions \eqref{def-kappa-delta}. 
It follows that 
\[ c_p(\s,\k,\d)  
\le \begin{cases}
1.014351 & \textrm{ if } \ p=2,\\
0.533948 & \textrm{ if } \ p\ge3,
    \end{cases}
\]
Note that since there is no primitive characters modulo $2$, there is no modulus $q$ with primitive characters of the form $2q'$ with $q'$ odd.
Thus $\nu_2\left(\frac{q}{q_k}\right) \ge 2$.
Since 
$\s\ge 0.9, 0.5 \le \d \le 0.75, 0.25 \le \k \le 0.5$, then
\[
 \frac{1-\kappa}{2c_p(\s,\k,\d)}\nu_p\left(\frac{q}{q_k}\right)  -  1 
\ge \begin{cases}
-0.507074 & \textrm{ if } \ p=2,\\
-0.531790 & \textrm{ if } \ p\ge 3.
    \end{cases}
\]
and
\[ C_p(\s,\k,\d) \ge 
a_0 + 
\begin{cases}
\min\left(0, \frac{1-\kappa}{c_2(\s_0,\k,\d)} - 1 \right)\sum_{k=2}^{n_0}  a_k   & \textrm{ if } \ p=2,\\
\min\left(0, \frac{1-\kappa}{2c_3(\s_0,\k,\d)} - 1 \right)\sum_{k=2}^{n_0} a_k  & \textrm{ if } \ p\ge 3.
\end{cases}
\]
We assumed that these quantities were positive, which achieves the proof.
\end{proof}
 \bibliographystyle{plain}
 \bibliography{Zeros-Dirichlet-Lfunctions-2017}
\end{document}